\documentclass[a4paper]{article}

\usepackage{amsmath}
\usepackage{amssymb}
\usepackage{mathrsfs}
\usepackage{amsthm}
\usepackage{tabularx}
\usepackage{longtable}
\usepackage{textcomp}
\usepackage{floatpag}  
\usepackage{enumerate} 

\newtheorem{thm}{Theorem}[section]
\newtheorem{defn}[thm]{Definition}
\newtheorem{lem}[thm]{Lemma}
\newtheorem{prop}[thm]{Proposition}

\newtheorem{rmk}[thm]{Remark}

\usepackage[top=1in, bottom=1.25in, left=1.25in, right=1.25in]{geometry}

\floatpagestyle{empty}

\newcommand{\ml}[2]{\begin{tabular}{@{} >{$}#1<{$} @{}} #2 \end{tabular}}

\title{Minimal $3$-generated Majorana algebras}
\date{}
\author{Andrey Mamontov\thanks{Sobolev Institute of Mathematics and Novosibirsk State University; mamontov@math.nsc.ru}, Alexey Staroletov\thanks{Sobolev Institute of Mathematics and Novosibirsk State University; staroletov@math.nsc.ru} 
and Madeleine Whybrow\thanks{Arbeitsgruppe Algebra, Geometrie und Computeralgebra, TU Kaiserslautern; mlw10@ic.ac.uk}}

\begin{document}

\maketitle

\begin{abstract}
Majorana theory was introduced by A.~A.~Ivanov as the axiomatization of certain properties of the $2A$-axes of the Griess algebra. Since its inception, Majorana theory has proved to be a remarkable tool with which to study objects related to the Griess algebra and the Monster simple group. We introduce the definition of a minimal $3$-generated Majorana algebra and begin the first steps towards classifying such algebras.

In particular, we give a complete classification of finite minimal $3$-generated $6$-transposition groups. We then use an algorithm developed in GAP by M.~Pfeiffer and M.~Whybrow, together with some additional computational tools, to give an almost complete description of all minimal $3$-generated Majorana algebras arising from this list of groups.
\end{abstract}

\section{Introduction}

The Monster group $\mathbb{M}$ is the largest of the twenty-six sporadic groups and was first constructed by Griess \cite{Griess82} as the automorphism group of the Griess algebra $V_{\mathbb{M}}$,
a $196884$-dimensional commutative non-associative real algebra. It is well known that $\mathbb{M}$ is a $6$-transposition group;
it is generated by a normal set of involutions such that the order of the product of any pair of these involutions is not greater than~6. Following the notation of \cite{Atlas85},
we refer to these involutions as the $2A$-involutions of the Monster.

Conway \cite{Conway85} showed that there exists a bijection $\psi$ between the $2A$-involutions of the Monster and certain non-trivial idempotents of the Griess algebra known as $2A$-axes.
This bijection leads to a deep connection between $6$-transposition subgroups of $\mathbb{M}$ and subalgebras of $V_{\mathbb{M}}$ generated by $2A$-axes.

For example, given any pair of $2A$-involutions of the Monster, their product lies in one of nine conjugacy classes in $\mathbb{M}$. Norton \cite{Norton96} classified all dihedral subalgebras of $V_{\mathbb{M}}$,
i.e. those generated by a pair of $2A$-axes. He showed that if $t_0$ and $t_1$ are $2A$-involutions of the Monster then the subalgebra generated by $\psi(t_0)$ and $\psi(t_1)$ has one of nine isomorphism types and that
this is determined by the conjugacy class of $t_0t_1$. This result provides a crucial tool with which to study the structure of the Griess algebra.

Later, Sakuma \cite{Sakuma07} reproved Norton's classification of the dihedral subalgebras of $V_{\mathbb{M}}$ in the more general context of Vertex Operator Algebras, or VOA's. Inspired by Sakuma's work, Ivanov \cite{Ivanov09}
axiomatized the properties of VOA's that were used in this result. This approach is known as Majorana theory and can also be thought of as a generalization of the relationship between the $2A$-involutions in the Monster
and the $2A$-axes in the Griess algebra.

The objects at the center of Majorana theory are real, commutative, non-associative algebras known as Majorana algebras. Each Majorana algebra is generated by certain idempotents known as Majorana axes which are in
bijection with involutions in the automorphism group of the algebra known as Majorana involutions.

In the seminal result of the theory, Ivanov et al. \cite{IPSS10} reproved Sakuma's theorem in the language of Majorana theory. In particular, they showed that a Majorana algebra generated by two axes must be isomorphic to one of the nine dihedral subalgebras of $V_{\mathbb{M}}$.

After this classification, it becomes natural to study Majorana algebras generated by three axes. In the following, we introduce the definition of a minimal
$3$-generated Majorana algebra and make the first steps towards classifying such algebras. In Section~2 we explain the concept of minimal $3$-generated $6$-transposition groups and classify all such finite groups.
In Section~3 we give some basic results from Majorana theory before presenting our work towards a classification of minimal $3$-generated algebras.

In particular, we have constructed $34$ minimal $3$-generated algebras, details of which are given in Table 7. Of the groups classified in Section~2, there are just three cases where we have not been able to classify all Majorana algebras that are admitted by a particular group. The outstanding cases are given in Table 8.

\section{Miminal 3-generated 6-transposition groups}
\label{sec:groups}

Let $G$ be a group. An element of $G$ is called an involution if its order divides $2$. Recall that a normal set of involutions $D$ in $G$ 
is called a set of 6-transpositions if $G = \langle D\rangle$ and, for each pair $a$ and $b$ from $D$, 
the order of $ab$ does not exceed 6. In this case, we say that $(G, D)$ is a 6-transposition group.

It is easy to see that if $(G, D)$ is a 6-transposition group generated by two elements
of $D$ then $G$ is a dihedral group of order not greater than $12$. However, 6-transposition groups
that can be generated by three elements of $D$ are not yet classified. Throughout this section
we deal with the following set of 3-generated 6-transposition groups.

\begin{defn}
Let $(G,D)$ be a 6-transposition group. 
We say that $(G,D)$ is \emph{minimal 3-generated} if it satisfies the following two properties:

\begin{enumerate}[(i)]
\item $G$ is generated by three elements of $D$;
\item $H\leq G$ and $H=\langle H\cap D \rangle$ then either $H=G$ or $H$ can be generated by two elements of~$D$.
\end{enumerate}
\end{defn}

Note that if $(G,D)$ is a minimal $3$-generated $6$-transposition group then $G$ need not be generated by three \emph{distinct} elements of $D$. For example, the trivial group and the dihedral groups $D_{2n}$ for $n \in \{1, 2, 3, 4, 5, 6\}$ obey this definition. Moreover, the class of minimal $3$-generated $6$-transposition groups is closed under homomorphic images. This will become especially important when we consider minimal $3$-generated Majorana algebras (see Section \ref{sec:constructions}). 

It is easy to see that every finite 6-transposition group $(G,D)$ generated by three distinct non-identity
elements of $D$ contains a minimal 3-generated group that is also generated by three distinct non-identity elements of~$D$.

\noindent\textbf{Notation.} The prime power $p^k$ denotes the elementary abelian group of order $p^k$, the number $n$ denotes the cyclic group of order $n$, 
$p^{1+2}$ denotes the extraspecial group of order $p^3$ which has no elements of order $p^2$. For groups $A$ and $B$, we write $A:B$ as well as $A \rtimes B$ to denote a split extension of $A$ by $B$, 
and $A \times B$ for their direct product.

In this section, we classify all finite minimal 3-generated 6-transposition groups. The main result is the following.
\begin{thm}\label{t:Min3GenGroups} Let $(G,D)$ be a non-trivial finite minimal 3-generated 6-transposition group. 
Then $G$ is isomorphic to one of the following.

\begin{enumerate}[(i)]
\item $D_{2n}$, where $n\in\{1,2,3,4,5,6,10\}$;
\item $2^3$, $S_4$, $GL(2,3)$, $A_5$;
\item $p^2:2\cong\langle x, y, z~|~x^p=y^p=[x,y]=z^2=x^zx=y^zy=1\rangle$, where $p\in\{3,5\}$;
\item $(2m \times 2):2\cong\langle x, y, z~|~x^{2m}=y^2=[x,y]=z^2=x^zx=y^zx^my=1\rangle$, \newline where 
$m=2$ and $D=y^G \cup z^G \cup (xz)^G$ or $D=y^G \cup z^G \cup (xz)^G \cup(x^2)^G$; \newline or $m=3$ and $D=y^G \cup z^G$ or $D=y^G \cup z^G \cup (x^3)^G$ ;
\item $p^{1+2}:2\cong\langle x, y, z~|~x^p=y^p=[x,y]^p=[x,[x,y]]=[y,[x,y]]=z^2=x^zx=y^zy=1\rangle$, where $p\in\{3,5\}$;
\item $(5^2:3):2\cong\langle x, y,z, w~|~x^5=y^5=z^3=[x,y]=w^2=x^wy^{-1}=y^zyx^{-1}=z^wz=1\rangle$.
\end{enumerate}
\end{thm}
\begin{rmk} In our proof we need the assumption that $G$ is a priori finite only in the cases when the product of any two distinct involutions from $D\setminus\{1\}$ has order 5, see Section~2.1.
\end{rmk}

\begin{rmk} The set $D \setminus \{1\}$ is unique if and only if all elements of order $2$ are conjugate in $G$. This is the case for items (iii), (v), (vi), for $G \cong D_{2n}$ with odd $n$ and $G=A_5$.
\end{rmk}

\subsection{Minimal 3-generated groups: preliminaries}

First we prove the following convenient criterion that helps to verify whether a given 6-transposition group
is minimal 3-generated or not.

\begin{lem}\label{l:criterion} Let $(G,D)$ be a 6-transposition group, where $G$ is generated by three elements of $D$. 
	Then $(G,D)$ is minimal 3-generated if and only if for every $a,b,c\in D$ the group $H=\langle a,b,c\rangle$ 
	either equals $G$ or can be generated by two elements of $D$.
\end{lem}
\begin{proof} If $(G,D)$ is minimal 3-generated and $a,b,c\in D$ then 
	the conclusion for $H=\langle a,b,c\rangle$ is true by the definition of minimal 3-generated groups.
	Now we prove the reverse implication. Let $H$ be a proper subgroup of $G$ 
	and $H=\langle H\cap D\rangle$. Take a generating set $S\subseteq D$ for $H$ of minimal size.
	Suppose that $|S|\geq3$ and take distinct $a$, $b$, and $c$ from $S$. Then $\langle a,b,c\rangle\leqslant H$
	and hence $\langle a,b,c\rangle\neq G$. By assumption, we know that $\langle a,b,c\rangle$ can be generated by two elements of $D$
	and so the size of $S$ is not minimal; a contradiction. Thus $|S|\leq 2$ and $G$ is minimal 3-generated.
\end{proof}

If $(G,D)$ is a 6-transposition group then we often write $G$ instead of $(G,D)$. In such cases, either $D$ is known or we only need the fact that a set $D$ exists. 

Each non-trivial minimal 3-generated group belongs to one of the following sets.

\begin{defn} Denote by $S(k,l,m)$ the set of minimal 3-generated 6-transposition groups $(G,D)$
that can be generated by three involutions $a$,$b$,$c\in D \setminus \{1\}$ with $|ab|=k$, $|bc|=l$, $|ac|=m$ and $1\leq k\leq l\leq m\leq 6$.
\end{defn}\label{sklm}

Note that a minimal 3-generated 6-transposition group can belong to different sets $S(k,l,m)$. 
For example, it is easy to see that $S_4\in S(2,3,3)$ and $S_4\in S(2,3,4)$.

We will describe minimal 3-generated groups using a simple observation that if $G\in S(k,l,m)$
then $G$ is a homomorphic image of the group $\langle a,b,c~|~R(k,l,m)\rangle$ where $$R(k,l,m)=\{a^2,b^2,c^2,(ab)^k,(bc)^l,(ac)^m\}.$$

Since dihedral groups are closely related with minimal 3-generated groups,
we need the following lemma.

\begin{lem}\label{l:dihedral} Let $(G,D)\in S(k,l,m)$ where $G$ is generated by two $D$-elements. Then $G\cong D_{2n}$
and one of the following claims holds:
\begin{enumerate}[(i)]
\item $n\in\{1,2,3,5\}$ and $(k,l,m)\in\{(1,n,n), (n,n,n)\}$;
\item $n=4$ and $(k,l,m)\in\{(1,4,4), (2,2,4), (2,4,4)\}$;
\item $n=6$ and $(k,l,m)\in\{(1,6,6), (2,2,3), (2,2,6), (2,3,6), (3,6,6)\}$.
\end{enumerate}
\end{lem}

\begin{proof}
It is clear that $n\leq 6$. Consider generators $a,b,c$ of $G$ as in Definition \ref{sklm}. We can assume they are distinct.  

Suppose $\langle a,c \rangle =G$. If $n$ is odd then $b \in a^G$ and $k=l=n$. If $n$ is even then either $b \in a^G \cup c^G$ and one of the possibilities holds, or $b \in Z(G)$ and $k=l=2$. 

Suppose now that $\langle a,c \rangle < G$. Then $n=6$ and $m=3$. It follows that $a$ and $c$ are conjugate. If $l=3$ then $b$ and $c$ are conjugate, but $D_{12}$ is not generated by one class of involutions. 
Therefore $k=l=2$, and item (iii) holds.
\end{proof}

Observe that the free Burnside group $B(2,5)$ of exponent $5$, generated by two elements $x$ and $y$, has an automorphism $\tau$ of order 2
such that $x^\tau=x^{-1}$ and $y^\tau=y^{-1}$. It is easy to see that the natural semidirect product $G=B(2,5)\rtimes\langle\tau\rangle$
is a 6-transposition group where $D$ is the set of elements of order 2. Moreover, the product of any two distinct involutions from $D$
has order 5 and $G$ is generated by three involutions $x\tau$, $y\tau$ and $\tau$. Since it is still unknown whether $B(2,5)$ is finite or not, 
the 6-transposition groups, where every $a,b \in D\setminus\{1\}$ satisfy $(ab)^5=1$, are of particular interest.
We introduce the following set of minimal 3-generated groups, which will be given special attention.
\begin{defn}\label{d:Kp} Let $p$ be an odd prime. Denote by $\mathcal{K}_p$ the class of groups $G$ such that
\begin{enumerate}[(i)]
\item $G$ includes a normal set of involutions $D$ such that $G$ is generated by three elements of $D$;

\item for every $x, y\in D \setminus \{1\}$ the order of $xy$ divides $p$;

\item if $H\leq G$ and $H=\langle H\cap D \rangle$ then either $H=G$ or $H$ can be generated by two elements of~$D$.
\end{enumerate}
\end{defn}

\subsection{Proof of Theorem~\ref{t:Min3GenGroups}: case $G$ is not a group from $\mathcal{K}_5$}
In this section, we prove Theorem~\ref{t:Min3GenGroups} for minimal 3-generated 6-transposition groups $(G,D)$
that do not belong to $\mathcal{K}_5$.

Throughout we suppose that $(G,D)$ is a non-trivial minimal 3-generated 6-transposition group.
By Lemma~\ref{l:dihedral}, we can assume that if $G$ is generated by three elements $a,b,c\in D$
then these three involutions are all distinct.

Our strategy is the following: first in Lemmas~\ref{l:2lm} and~\ref{l:255} we consider groups $G$ generated by $a,b,c\in D$ with commuting $a$ and $b$. 
Then in Lemma~\ref{l:5} we consider the case, when there exist $e,f\in D$ such that $|ef|=5$. Finally, we proceed with $G\in S(k,l,m)$, where $k\in\{3,4,6\}$. This order of cases allows us
to use a reduction.

In each of these cases, the fact that $(G,D)$ is a minimal $3$-generated $6$-transposition group forces the elements of $G$ to obey relations in addition to those contained in $R(k,l,m)$. These new relations either force certain 3-generated subgroups of $G$ to be dihedral or are of the form $(xy)^r$ for $x,y \in D$ and $r \in \{1, \ldots, 6\}$. Using an algorithm for coset enumeration in GAP \cite{GAP16}, we verify that these sets of relations define finite groups. 
At this stage, we find all homomorphic images of these finite groups and use Lemma~\ref{l:criterion} to verify which of these are minimal $3$-generated.

\begin{lem}\label{l:2lm}
If $G\in S(2,l,m)$, where $(l,m)\neq(5,5)$, then $G$ is isomorphic to one of the following: $2^3$, $S_4$, $A_5$ or $D_{2n}$ with $n\in\{2,4,6,10\}$.
\end{lem}
\begin{proof}

Let $l=2$. Then $b \in Z(G)$ and clearly relations $R(2,2,m)$ define the group $2\times D_{2m}$.
Groups $2\times D_8$ and $2\times D_{12}$ are not minimal 3-generated, 
since in these cases $\langle b, a, a^c \rangle \cong 2^3$ and $\langle b, a, c^{ac} \rangle \cong 2^3$ are proper non-dihedral subgroups, respectively.
So $G$ is one of the following: $D_4$, $2^3$, $2 \times D_6 \cong D_{12}$, $D_8$ or $2 \times D_{10} \cong D_{20}$.

Suppose that $l = 3$ or $4$ and $m \neq 6$. Where necessary, we use the fact that $G$ is a $6$-transposition group to derive relations in addition to those coming from $R(2, l, m)$. These additional relations are shown in the second column of Table \ref{t:main}. Using coset enumeration, we see that these extended sets of relations define the groups given in the third column of this table. We then calculate the homomorphic images of these groups and verify which are minimal 3-generated with the aid of Lemma~\ref{l:criterion}. In these cases, we infer that $G\in\{D_8, S_4, A_5\}$.

Suppose that $l=3$ and $m=6$. Since $ab=ba$, we have $|bc^a|=|b^ac|=|bc|=3$. Also we have $|cc^a|=|(ca)^2|=3$. 
If $\langle b,c^a,c \rangle$ is generated by two $D$-elements then Lemma~\ref{l:dihedral} implies that $b \in \langle c^a, c \rangle \subset \langle a,c \rangle$ and hence $G \cong D_{12}$. 
So we can assume that $G = \langle b,c^a,c \rangle \in S(3,3,3)$. However, the minimal $3$-generated groups contained in $S(3,3,3)$ are listed in Table~\ref{t:main}. As none of these contain an element of order $6$ and a pair of distinct commuting involutions, they cannot also lie in $S(2,3,6)$ and so this case cannot occur. 

It remains to consider the cases $(l,m)\in\{(4,6),(5,6),(6,6)\}$. Assume that $a,b,c\in D$ such that $|ab|=2$, $|bc|=l$, where $l\in\{4,5,6\}$, and $|ac|=6$.
Put $w=a^{ca}$. Then $|aw|=|aacaca|=|(ca)^2|=3$. Consider the group $H=\langle a,b,w \rangle$. 
If $|wb| \in \{3,4,5 \}$ then $H\in S(2,3,|wb|)$, in particular, $H$ is not dihedral by Lemma~\ref{l:dihedral}. So $G=H$ and these cases are considered above.
Therefore, either $|wb|=2$ or $|wb|=6$. If $|wb|=2$ then using coset enumeration, we see that $G$ is a homomorphic image of $D_8\times S_3$, $D_{20}$ or $2\times S_3\times S_3$.
As above, we verify that there are no groups belonging to $S(2,l,6)$ among these images.

Finally, let $|wb|=6$. We have shown that $S(2,3,6)=\{D_{12}\}$. So we can assume that $H$ is generated by two $D$-elements
and hence $H=\langle b,w \rangle$. Since $|ab|=2$ and $|aw|=3$, we conclude that $a=w^{bw}$. After adding this relation to $R(2,l,6)$, 
we see that if $l=4,5,6$ then $G$ is a homomorphic image of $2\times D_8$, $D_{20}$ or $2\times S_3\times S_3$, respectively. The first two of these groups have no elements of order $6$.
Homomorphic images of $2\times S_3\times S_3$ either contain $D$-generated subgroup isomorphic to $2^3$ or do not belong to $S(2,6,6)$. So this case is not possible.
 
\end{proof}

\begin{lem}\label{l:255} If $G\in S(2,5,5)$ then $G\cong A_5$.

\end{lem}

\begin{proof}
Let $G=\langle a,b,c\rangle$, $|ab|=2$ and $|bc|=|ac|=5$.
Since $(ac)^5=(bc)^5=1$, we have $c^{acac}=a$ and $c^{bcbc}=b$. The order of $ab$ is two, so $a$ and $b$ commute. Hence, their conjugates $c=b^{cbcb}$ and $a^{cbcb}$ commute.  
It follows that $a^{cbcbacac}$ and $c^{acac}=a$ commute.
Let $w=a^{cbcbacac}$. If $w=a$ or $w=b$ then using coset enumeration in GAP, we infer that $G$ is a homomorphic image of either $D_{10}$ or a group of form $((2\times Q_8):2):D_{10}$.
The first case contradicts to Lemma~\ref{l:dihedral}. In the second case put $x=a^{cbcb}$. Then verifying images, we see that $\langle a^x, a^c, b^c \rangle$ is isomorphic either to $2\times D_8$ or $2^3$,
so there are no minimal 3-generated groups in this case.

Therefore, we can assume that $a$, $b$ and $w$ are distinct elements of $D$. Let $H=\langle b,a,w\rangle$. We know that $H\in S(2,2,|wb|)$ and hence if $G=H$ then $G\cong A_5$ or $G\cong D_{20}$ by Lemma~\ref{l:2lm}.
It is easy to see that $D_{20}\not\in S(2,5,5)$, so $H$ is generated by two $D$-elements.

Since $a$ commutes with $w$ and $b$, we conclude that $a\in Z(H)$ and hence $H\in\{D_{4}, D_8, D_{12}\}$. Then $H$
contains exactly three conjugacy classes of elements of order 2, and since $H$ is generated by two $D$-elements, all three classes are included in $D$.
Therefore, $ab\in D$. Assume that $|abc|\neq5$. Then the group $\langle a,ab, c \rangle$ is not dihedral by Lemma~\ref{l:dihedral}
and hence $\langle a,ab, c \rangle=G$. Therefore $G$ is isomorphic to one of the groups listed in the conclusion of Lemma \ref{l:2lm}. Since $|G|$ is divisible by 5 and $G\in S(2,5,5)$, we have $G\cong A_5$.
So $|abc|=5$. Considering the group $\langle a,b^c,b \rangle$ and reasoning in the similar way, we infer that $|ab^c|=5$. Using coset enumeration, 
we see that the set of relations $R(2,5,5)\cup\{(abc)^5, (ab^c)^5\}$ defines the trivial group.

\end{proof}

\begin{lem}\label{l:5} Suppose that there exist $a,c \in D$ with $|ac|=5$. Then $G$ is isomorphic to one of the following: $D_{10}$, $D_{20}$, $A_5$, or $(5^2:3):2$.
\begin{proof}
Let $(G,D)$ be a counterexample. 
If there exists $d \in D\setminus\{1,a\}$ such that $[a,d]=1$ then Lemma~\ref{l:dihedral} implies that $\langle a,c,d \rangle$ is non-dihedral and therefore $G=\langle a,c,d \rangle$;
a contradiction with Lemmas \ref{l:2lm} and \ref{l:255}.
So we can assume that for all $d \in D\setminus\{1\}$ the order of $ad$ is odd, in particular, $a$ and $d$ are conjugate in $\langle a,d\rangle$. 
Therefore, the order of the product of any two distinct elements in $D\setminus\{1\}$ is odd.

If for every $b \in D$ a subgroup $\langle a,b,c \rangle$ can be generated by two $D$-elements then $b \in \langle a,c \rangle$, so $G \cong D_{10}$, a contradiction.
Thus for some $b \in D$ we have $G= \langle a,b,c \rangle$. 

First, assume that $G \in S(3,3,5)$ and choose $b\in D$ with $|ab|=|bc|=3$. Then $G$ is a homomorphic image of groups $G(i,j)=\langle a,b,c ~|~R(3,3,5),(ba^c)^i,(ba^{ca})^j \rangle$, where $i,j \in \{3,5\}$. 
Using coset enumeration, we see that $G(3,5)\cong G(5,3)\cong 2$, and $G(5,5) \cong 2 \times PSU_3(4)$. If $G\cong 2 \times PSU_3(4) $ or $G\cong PSU_3(4)$ then
calculations show that there exists an element of order 2 in $D\setminus\{a\}$ commuting with $a$; a contradiction.
Finally, we see that $|G(3,3)|=5^2 \cdot 3 \cdot 2$. In this case, put $x=ac$, $y=x^b$, $z=ab$ and $w=b$. 
Then defining relations of case (vi) of Theorem \ref{t:Min3GenGroups} hold and so $G(3,3) \cong (5^2:3):2$; a contradiction.

Assume further that $|ab|=3$ and $|bc|=5$. Then $a=b^{ab}$ and $b^{(cb)^2}=c=a^{(ca)^2}=b^{ab(ca)^2}$. 
It follows that $x=(ab) (ca)^2 \cdot (bc)^2 \in C_G(b)$. Conjugating $(ab)^3=1$, we obtain $(a^xb)^3=1$. So if $G=\langle a,b,a^x \rangle$ then $G \in S(3,3,5)$ or $G\in S(3,3,3)$; a contradiction.
Therefore, we can assume that $\langle a,b,a^x \rangle$ is generated by two $D$-elements. Then $\langle a, b \rangle\cong S_3$ and $a^x \in\{a,b,a^b\}$. If $a^x =b$ then $a=xbx^{-1}=b$; a contradiction.
Therefore either $a^x=a$ or $a^x=a^b$.

Let $b_1=b^c$, $b_2=b^{cb}$ and $b_3=b^{cbc}$ be distinct involutions of $\langle b,c \rangle$. Note that $|bb_j|=5$ and $|ab|=3$, so $H=\langle b_j,a,b \rangle $ is not generated by two $D$-elements.
Therefore, if $|ab_j|=3$ for some $j$ then $G=H \in S(3,3,5)$; a contradiction. Hence, we can assume that $|ab_j|=5$.

Thus $G$ is a homomorphic image of groups $G(i)=\langle a,b,c~|~R(3,5,5),(ab_1)^5,(ab_2)^5,(ab_3)^5,a^xa^{b^i} \rangle$ where $i=0$ in the case $a^x=a$, or $i=1$ in the case $a^x=a^b$. 
A group $G(0) \cong 3^4 : D_{10}$ is not minimal 3-generated, because $\langle a,b,a^{cabc} \rangle \in S(3,3,3)$ is isomorphic to $3^2 : 2$, and the orders of other homomorphic images of $G(0)$ 
are not divisible by 3. At the same time, $G(1)\cong D_{10}$ and again $|G(1)|$ is not divisible by 3; a contradiction. 
\end{proof}
\end{lem}

\begin{lem}\label{l:246} If $G\in S(3,l,m)$ and for every $x,y\in D$ it is true that $|xy|\neq5$ then $G$ is isomorphic to one of $S_3$, $D_{12}$, $3^2:2$, $S_4$, $GL(2,3)$ or $3^{1+2}:2$.
 \end{lem}
\begin{proof} Let $a,b,c\in D$ such that $|ab|=3$, $|bc|=l$, $|ac|=m$ and $G=\langle a,b,c \rangle$. First we consider the case $l=m=3$.
We add the relation $(ab^c)^j=1$, where $j\in\{4,6\}$, to $R(3,3,3)$. Using coset enumeration and Lemma~\ref{l:criterion}, we see that $G$ is isomorphic to one of $S_3$, $S_4$, $3^2:2$, $3^{1+2}:2$ (see Table~\ref{t:main}).
So we can assume that $m=4$ or $6$.

The strategy is to consider all possibilities of $(l,m)$, and in each case choose a word $w$ such that $\langle a,b,w \rangle$ has already been considered.

Suppose that $(l,m)\in\{(3,4),(4,4)\}$. Put $w=a^c$. Since $|ac|=4$, we have $|aw|=2$. Consider the group $H=\langle w,a,b \rangle$.
If $|bw|=3,4$ then $H\in S(2,3,3)$ or $H\in S(2,3,4)$. By Lemmas~\ref{l:dihedral} and \ref{l:2lm}, we conclude that $G=H\cong S_4$. Therefore $|bw|=2$ or $|bw|=6$.
In the first case, we see that if $l=3,4$ then $G$ is a homomorphic image of $2$ or $(S_3\times S_3):2$, respectively. 
Depending on the set $D$, the group $(S_3\times S_3):2$ contains either $D_{12}$ or $S_3\times S_3$ which are 3-generated and cannot be generated by two $D$-elements. 

Therefore, we can assume that $|wb|=6$. Then $H\in S(2,3,6)$ and hence either $H=G$ or $H=D_{12}$. In the first case, Lemma~\ref{l:2lm} implies  that $G=D_{12}$,
which is not possible, since $D_{12}$ has no elements of order 4. So $H=D_{12}$ and hence $H=\langle w,b\rangle$. Since $|ab|=3$ and $|aw|=6$, we have $a=b^{wb}$. 
Using this relation, we see that $G$ is a homomorphic image of $GL(2,3)$ or $(2^2 \times A_4): 2$. 
In the first case, we obtain $G=GL(2,3)$, since for other homomorphic images $|ac|<4$. In the second case, if $G=(2^2 \times A_4): 2$ then
there exist three elements in $a^G$ that generate a subgroup isomorphic to $2^3$, so $G$ is not minimal 3-generated. 
Further, we see that there is only one other homomorphic image such that $ab$, $bc$ and $ac$ have same orders as in the initial group. This image is isomorphic to $(6\times 2):2$.
To see that this group is isomorphic to the corresponding group in Theorem~\ref{t:Min3GenGroups}, one can take $x=bw$, $y=c$, $z=a$ and verify that all required relations hold true for $x$, $y$ and $z$.

Suppose that $l=3$ and $m=6$. Put $w=c^{ac}$ and $H=\langle a,b,w \rangle$. As above, we can assume that either $|bw|=2$ or $|bw|=6$. 
In the first case, $G$ is a homomorphic image of $S_3$; a contradiction. So $|bw|=6$. Then $H\in S(2,3,6)$ and $H=D_{12}$ by Lemma~\ref{l:2lm}.
Therefore $a=b^{wb}$. We use this equality and relation $(ab^c)^j=1$ where $j=4$ or $j=6$. Applying coset enumeration, we see that 
if $j=6$ then $G$ is finite and $(ac)^3=1$; a contradiction. 
Similarly, if $j=4$ then $G$ is a homomorphic image of $(((4\times 2):4):3):2$. For the latter group, there exist three elements in $a^G$ that generate a subgroup isomorphic to $2^3$.
Calculations show that there is only one other homomorphic image such that $ab$, $bc$ and $ac$ have same orders as in the initial group. This image is isomorphic to $GL(2,3)$.

Finally, let $(l,m)\in\{(4,6),(6,6)\}$. Put $w=a^c$ and $H=\langle a,b,w \rangle$. Then $|aw|=|acac|=3$. If $H=G$ then $G\in S(2,3,3)$ or $G\in S(3,3,i)$.
So we can assume that $H$ is generated by two $D$-elements. It follows from Lemma~\ref{l:dihedral} that $|bw|=3$ and $H=\langle a,b\rangle$. Then $w=a^b$. Using this relation,
we see that $G$ is a homomorphic image of $D_{12}$.
\end{proof}

\begin{lem}\label{l:even}
If $G\in S(k,l,m)$, where $k,l,m \in \{4,6\}$, then either $G \in S(u,v,w)$, with $u \in \{2,3\}$, or $G \cong (4 \times 2 ):2$.
\end{lem}

\begin{proof}
First, assume that $m=6$ and consider $x=c^{ac}$. Then $|ax|=2$ and $|cx|=3$. We can assume that both groups $\langle a,x,b \rangle$ and $\langle b,x,c \rangle$ are generated by two $D$-elements.
If $l=4$ then this is not possible by Lemma~\ref{l:dihedral}. Considering $a^{ca}$ instead of $x$, we similarly see that if  $G \in S(4,6,6)$ then $G\in S(u,v,w)$ where $u \in \{2,3\}$.
Assume now that $k=6$. By Lemma~\ref{l:dihedral}, we have $\langle a,x,b \rangle$ and $\langle b,x,c \rangle$ are isomorphic to $D_{12}$. Then $|bx|=2$ and $x=(ab)^3=c^{bc}$. 
Hence, $G$ is a homomorphic image of $\langle a,b,c ~|~ R(6,6,6),(bx)^2,x(ab)^3,xc^{bc}\rangle \cong D_{12}$ and $|bc|$ divides 2; a contradiction.

Finally, we assume that $k=l=m=4$. Let $x=a^b$. Then $|ax|=2$, and without loss of generality suppose that $\langle a,x,c \rangle=\langle a,c \rangle$ is dihedral. 
Lemma~\ref{l:dihedral} implies that either $|xc|=4$ and $x=a^c$, or $|xc|=2$ and $x=(ac)^2$. So $G$ is a homomorphic image of 
$G(i)=\langle a,b,c ~|~ R(4,4,4), xa^ia^c\rangle$ where $x=a^b$ and $i \in \{0,1\}$. 
Note that $G(1) \cong D_8$ and $a=1$ in this case; a contradiction. Moreover, $G(0) \cong (4 \times 2^2):2$ and this group contains a subgroup isomorphic to $2^3$ which is generated
by three $D$-elements. Inspecting other homomorphic images, we find for only one of them it is true that $|ab|=|bc|=|ac|=4$. This image is isomorphic to $(4 \times 2 ):2$.
\end{proof}

\begin{table}
\begin{center}
\begin{tabular}{|c|c|c|c|} \hline
$(k,l,m)$ & extra relations & $G$ &  \begin{tabular}{c} quotients of $G$ \\ belonging to $S(k,l,m)$ \end{tabular}   \\ \hline
$(2,3,3)$ & no & $S_4$ & $S_4$ \\ 
$(2,3,4)$ & no & $2\times S_4$ & $S_4$ \\ 
$(2,3,5)$ & no & $2\times A_5$ & $A_5$ \\ 
$(2,4,4)$ & $(ab^c)^6$ & $(2^2 \times S_3^2):2$ & $D_8$ \\ 
          & $(ab^c)^5$ & $(D_{10}\times D_{10}) : 2$ & $D_8$  \\ 
          & $(ab^c)^4$ & $(D_8 \times D_8) : 2$ & $D_8$  \\ 
$(2,4,5)$ & $(ab^c)^5$ & $A_6 : 2$ & --- \\ 
	  & $(ab^c)^4$ & $2 \times ((2^4 : 5) : 2)$ & --- \\ 
	  & $(ab^c)^6$, $(ca^{b^ca})^6$ &  $((3^4 : A_5) : 2) : 2$ & --- \\ 
	  & $(ab^c)^6$, $(ca^{b^ca})^5$ &  $D_{20}$ & --- \\ 
	  & $(ab^c)^6$, $(ca^{b^ca})^4$ &  $2 \times ((2^4 : A_5) : 2)$ & --- \\ 
$(3,3,3)$ & $(ab^c)^6$  & $((6\times 6) : 3) : 2$  & $S_3$, $3^2:2$, $S_4$, $3^{1+2} : 2$ \\ 
	  & $(ab^c)^5$  & $(5^2: 3) : 2$  & $(5^2 : 3) :2$, $S_3$ \\ 
	  & $(ab^c)^4$  & $((4\times4) :3) : 2$  & $S_3$, $S_4$ \\ \hline
\end{tabular}
\caption{3-generated groups and their minimal 3-generated quotients}
\label{t:main}
\end{center}
\end{table}

\subsection{Classification of finite groups in $\mathcal{K}_p$}

By results from the previous section, to prove Theorem~\ref{t:Min3GenGroups} it remains to
show that finite groups from $\mathcal{K}_5$ are listed in the theorem. In this section,
we classify all finite groups in $\mathcal{K}_p$ for arbitrary odd prime $p$.

For a finite group $H$, denote by $O(H)$ the largest normal subgroup of odd order in $H$
and by $Z^{\ast}(H)$ the inverse image in $H$ of the center of $H/O(H)$.
We need the following corollary from the famous $Z^\ast$-theorem.

\begin{lem}\label{Z-thm}\cite[Corollary 1]{Glaub66}
Let $S$ be a Sylow 2-subgroup of a finite group $G$. Suppose $x\in S$. 
A necessary and sufficient condition for $x\not\in Z^{\ast}(G)$ is that there exists $y\in C_S(x)$ such that $y$ is conjugate to $x$ in $G$ and $y\neq x$. 	
\end{lem}

Throughout this section we assume that $G$ is finite and $G\in\mathcal{K}_p$. 
By Definition~\ref{d:Kp}, there exists a normal set of involutions $D$ in $G$
and $a,b,c\in D$ such that $G=\langle a,b,c \rangle$. If $G \neq 1$ then we assume that $a\neq 1$.
The strategy is to prove that $O(G)$ is $p$-group, describe its structure, and use that to classify groups from $\mathcal{K}_p$.

\begin{lem}\label{l:s2:lem1} The following holds.
\begin{enumerate}[(i)]
 \item $G=O(G)\rtimes\langle{a}\rangle$, in particular, all non-identity involutions in $G$ are conjugate. 
Moreover, $O(G)=\langle ab,ac\rangle$;
\item If $N$ is a proper $\langle{a}\rangle$-invariant normal subgroup in $O(G)$ then $G/N\in~\mathcal{K}_p$. 
\end{enumerate}

\end{lem}
\begin{proof}
Let $T\in Syl_2(G)$ such that $a\in T$. If $g\in G$ and $a^g\neq a$ then $|aa^g|=p$, so $a^g\not\in T$. By Lemma~\ref{Z-thm},
we obtain $a\in Z^\ast(G)$. Similarly, we see that $b,c\in Z^\ast(G)$ and hence $Z^\ast(G)=G$. So, $|G/O(G)|$ divides 8 and $ab,bc,ac\in O(G)$.
Therefore, $|G/O(G)|=2$ and $G=O(G)\rtimes\langle{a}\rangle$.  Now we prove that $O(G)$ is generated by $ab$ and $ac$.
Observe that $bc=baac$, so $bc\in\langle ab, ac\rangle$. Let $x\in O(G)$ and
$x=i_1i_2...i_k$ where $i_j\in\{a,b,c\}$ for all $j$ and $i_j\neq i_{j+1}$ for $j<k$. Clearly $k>1$ and $i_1i_2\in\langle ab, ac\rangle$.
So $i_3i_4...i_k\in O(G)$. Using induction, we obtain  $x\in\langle ab,ac\rangle$ and hence $O(G)=\langle ab,ac\rangle$.

Now we prove $(ii)$. If $a\in G$ and $A\leq G$ then we denote the images of $a$ and $A$ under the canonical homomorphism from $G$ to $G/N$
by $\overline{a}$ and $\overline{A}$, respectively.
Let $\overline{t}$ be an involution in $G/N$. Then we can assume that $t$ is an involution. 
Indeed, $t^2\in O(G)$ and hence $|t|=2(2m+1)$ for $m\in\mathbb{Z}$,
so $|t^{2m+1}|=2$ and $\overline{t^{2m+1}}=\overline{t}$. Clearly $\overline{G}$ is generated by three involutions $\overline{a}$, $\overline{b}$ and $\overline{c}$.
As a normal set of involutions for $\overline{G}$, we can take $\overline{a}^{\overline{G}}$.
Let $\overline{t_1}\neq\overline{t_2}\in\overline{a}^{\overline{G}}$.
If $|\overline{t_1}|=|\overline{t_2}|=2$ and $t_1\neq t_2$ then $\overline{t_1\cdot t_2}=\overline{t_1}\cdot\overline{t_2}$, 
so $|\overline{t_1}\cdot\overline{t_2}|=p$. 
Let finally $\overline{t_1},\overline{t_2},\overline{t_3}\in\overline{a}^{\overline{G}}$ and  $\langle\overline{t_1},\overline{t_2},\overline{t_3}\rangle\neq\overline{G}$. 
Then, obviously, $\langle t_1,t_2,t_3\rangle\neq G$. So  $\langle t_1,t_2,t_3\rangle$ is a homomorphic image of $D_{2p}$ and hence $\langle\overline{t_1},\overline{t_2},\overline{t_3}\rangle$
is a homomorphic image of $D_{2p}$. Thus $\overline{G}\in\mathcal{K}_p$.
\end{proof}

\begin{lem}\label{l:s2:lem2} $O(G)$ is a $p$-group.
\end{lem}
\begin{proof}
We proceed by induction on $|G|$ and assume that $|G|>2$. Then $p$ divides $|O(G)|$.
Let $N=O_{p'}(O(G))\neq1$. Lemma~\ref{l:s2:lem1} implies that $G/N\in\mathcal{K}_p$ and hence $O(G)/N$ is a $p$-group by induction. 
Therefore, $|G:N|$ and $|N|$ are coprime, and the Schur--Zassenhaus theorem yields there exists a subgroup $T\leq G$ isomorphic to $G/N$. 
Since $T\in\mathcal{K}_p$ and $p$ divides $|T|$, we infer that $T\cong D_{2p}$.
Since $|N|$ is not divisible by $p$ and $G\in\mathcal{K}_p$, we have $a$ is the unique element of order 2 in $N\langle{a}\rangle$. 
So $N\leq C_G(a)$. Therefore, $|G/C_G(a)|$ divides $p$ and hence $T$ contains all involutions of $G$; a contradiction, since $G$ is generated by involutions.

So $N=O_{p}(O(G))\neq1$. Suppose that $N<O(G)$. Then $G/N\in\mathcal{K}_p$ and hence $O(G)/N$ is a $p$-group by the inductive hypothesis; 
a contradiction. 
\end{proof}

\begin{lem}\label{l:s2:lem3}
Let $x=ab$ and $y=ac$. Then the following claims hold.
\begin{enumerate}[(i)]
\item  If $|O(G)|=1$ then either $G=1$ or $G=Z_2$;	
\item  if $|O(G)|=p$ then $G=D_{2p}$;
\item if $|O(G)|=p^2$ then $G=(\langle x\rangle\times\langle y \rangle)\rtimes\langle{a}\rangle$, 
where $x^a=x^{-1}$, $y^a=y^{-1}$;
\item if $|O(G)|=p^3$ then $O(G)\cong p^{1+2}$ and
$G=\langle x,y \rangle \rtimes\langle{a}\rangle$, where $[x,y]^p=1$ and $[x,y]\in Z(G)$.
\end{enumerate}
\end{lem}
\begin{proof} If $|O(G)|\leq p$ then the assertions are clear. Let $|O(G)|=p^2$. Then we have $O(G)=\langle x,y\rangle$ by Lemma~\ref{l:s2:lem1}.
Since $|x|=|y|=p$, we infer that $O(G)$ is elementary abelian. Note that $(ax)^2=b^2=1$, so $x^a=x^{-1}$. Similarly $y^a=y^{-1}$
as required.

Finally, let $|O(G)|=p^3$. Since $O(G)=\langle x,y\rangle$, we have $O(G)$ is non-abelian.  Since $O(G)$ is generated by two elements of order $p$, 
it is isomorphic to $p^{1+2}$. Then $[x,y]\neq1$ and hence $\langle[x,y]\rangle=Z(O(G))=O(G)'$. 
Note that $x^a=x^{-1}$ and $y^a=y^{-1}$. So $(x^{-1}y^{-1}xy)^a=xyx^{-1}y^{-1}$.
Since $[x,y]\in Z(O(G))$, we see that $xy=(yx)(x^{-1}y^{-1}xy)=(x^{-1}y^{-1}xy)(yx)$. Therefore,  $xyx^{-1}y^{-1}=x^{-1}y^{-1}xy$ and hence $[x,y]^a=[x,y]$ as required.
\end{proof}

The following proposition completes the classification of finite groups in $\mathcal{K}_p$.

\begin{prop} Every group from $\mathcal{K}_p$ is isomorphic to a group described in Lemma~\ref{l:s2:lem3}.
\end{prop}
\begin{proof}
First, it is easy to verify that all groups from Lemma~\ref{l:s2:lem3} are in $\mathcal{K}_p$.
Suppose that there exists $G\in\mathcal{K}_p$ such that $|G|\geq 2p^3$. We show that $|G|=2p^3.$
Let $P=O(G)=\langle x,y \rangle$ be the Sylow $p$-subgroup of $G$. Since $P$ is generated by two elements of order $p$
and has order greater than $p^2$, we conclude that $P$ is non-abelian, in particular, it is not cyclic.
So the minimal number of generators of $P$ is $2$ and hence the Frattini subgroup $\Phi(p)$ has index $p^2$ in $P$.

Now we prove that $a$ centralizes $[x,y]$ and $|[x,y]|=p$. Since $\Phi(P)$ is characteristic in $P$, it is normal in $G$. 
Lemma~\ref{l:s2:lem1} implies that $G/\Phi(P)$ is isomorphic
to the group from Lemma~\ref{l:s2:lem3}(iii). Therefore, there exists a normal subgroup $\overline{N}$ in $G/\Phi(P)$ of order $p$.
Let $N$ be the full preimage of $\overline{N}$ in $G$ and $T=N\langle a\rangle$. Suppose that $T$ has more than $p$ distinct 
conjugates of $a$. Then we can take three distinct involutions $a$, $i$ and $j$ such that $j\not\in\langle a,i \rangle$.
Then $|\langle a,i,j \rangle|>2p$ and $|\langle a,i,j \rangle|\leq|T|<|G|$; a contradiction with $G\in\mathcal{K}_p$.
So $T$ has at most $p$ involutions and hence $|T:C_T(a)|\leq p$. So $|C_T(a)|\geq |\Phi(P)|$. Since there are no elements of order $2p$
in $G/\Phi(P)$, we conclude that $C_T(a)=\Phi(P)$. Since $P/\Phi(P)$ is abelian, we have $P'\leqslant \Phi(P)$, in particular, $P'\leqslant C_G(a)$.
Then $[x,y]^a=[x,y]$ and $[x,y^{-1}]^a=[x,y^{-1}]$. 

We know that $(x^{-1}y^{-1}xy)^a=(x^{-1})^a(y^{-1})^ax^ay^a=xyx^{-1}y^{-1}$.
Therefore $(xy)(x^{-1}y^{-1})=(x^{-1}y^{-1})(xy)$. Since $xy=xaay$ and $(ax)^2=(ay)^2=1$, we have $|xy|=p$.
Similarly, $|x^{-1}y^{-1}|=p$ and hence $[x,y]^p=(xy)^p(x^{-1}y^{-1})^p=1$ as claimed.
Considering the equality $[y^{-1},x]^a=[y^{-1},x]$, we infer that $(yx^{-1})(y^{-1}x)=(y^{-1}x)(yx^{-1})$.

Now we prove that $x[x,y]=x^{2n}y^{-1}xyx^{-2n}$
for every integer $n\geq0$. For $n=0$, we see that $x[x,y]=xx^{-1}y^{-1}xy=x^0y^{-1}xyx^0$.
Suppose that $x[x,y]=x^{2n}y^{-1}xyx^{-2n}$. Then 
\begin{multline*}
x[x,y]=x^{2n}y^{-1}xyx^{-2n}=x^{2n+1}(x^{-1}y^{-1})(xy)x^{-2n}=x^{2n+1}(xy)(x^{-1}y^{-1})x^{-2n}=\\
=x^{2n+2}yx^{-1}y^{-1}x^{-2n}=x^{2n+2}(yx^{-1})(y^{-1}x)x^{-2n-1}=x^{2n+2}(y^{-1}x)(yx^{-1})x^{-2n-1}=\\
=x^{2n+2}y^{-1}xyx^{-2n-2} \text{ as required.}
\end{multline*}
Put $n=p-1$, we obtain $x[x,y]=x^{p-1}y^{-1}xyx^{-(p-1)}=x^{-1}y^{-1}xyx=[x,y]x$.
So $[x,y]$ and $x$ commute. Similarly $[x,y]$ and $y$ commute. Therefore, $x^p=y^p=[x,y]^p=1$ and $[x,y]\in Z(P)$. So $P$ is a homomorphic
image of $p^{1+2}$, in particular, $|P|\leq p^3$. Thus $|P|=p^3$ and the proposition is proved.

\end{proof}

\section{Majorana theory}

Majorana theory was introduced by Ivanov in \cite{Ivanov09} as an axiomatisation of certain properties of the $2A$-axes of the Griess algebra. In this section, we introduce the definition of a minimal $3$-generated Majorana algebra. We then use our results from Section~\ref{sec:groups} to give the first steps towards a classification of minimal $3$-generated algebras.
\subsection{Definitions and background}

Let $V$ be a real vector space equipped with a positive definite symmetric bilinear form $( \, , \, )$ and a bilinear commutative non-associative algebra product $\cdot$ such that
\begin{description}
\item[M1] $( \, , \, )$  associates with $\cdot$ in the sense that
\[
(u , v \cdot w ) = (u \cdot v, w )
\]
for all $u,v,w \in V$;

\item[M2] the \emph{Norton inequality} holds so that
\[
(u \cdot u, v \cdot v) \geq (u \cdot v, u \cdot v)
\]
for all $u,v \in V$.
\end{description}
Let $A$ be a subset of $V \backslash \{0\}$ and suppose that for every $a \in A$ the following conditions M3 to M7 hold:
\begin{description}
\item[M3] $(a,a) = 1$ and $a \cdot a = a$, so that the elements of $A$ are idempotents of length 1;

\item[M4] $V = V_1^{(a)} \oplus V_0^{(a)} \oplus V_{\frac{1}{2^2}}^{(a)} \oplus V_{\frac{1}{2^5}}^{(a)}$ where $V_{\mu}^{(a)} = \{ v \mid v \in V, \, a \cdot v = \mu v\}$ is the set of $\mu$-eigenvectors of the adjoint action of $a$ on $V$;

\item[M5] $V_1^{(a)} = \{ \lambda a \mid \lambda \in \mathbb{R} \}$;

\item[M6] the linear transformation $\tau(a)$ of $V$ defined via

\[
\tau(a): u \mapsto (-1)^{2^5 \mu}u
\]

for $u \in V_{\mu}^{(a)}$ with $\mu = 1,0,\frac{1}{2^2}, \frac{1}{2^5}$, preserves the algebra product (i.e. $u^{\tau(a)} \cdot v^{\tau(a)} = (u \cdot v)^{\tau(a)}$ for all $u, v \in V$);

\item[M7] if $V_+^{(a)}$ is the centraliser of $\tau(a)$ in $V$, so that $V_+^{(a)} = V_1^{(a)} \oplus V_0^{(a)} \oplus V_{\frac{1}{2^2}}^{(a)}$, then the linear transformation $\sigma(a)$ of $V_+^{(a)}$ defined via
\[
\sigma(a) : u \mapsto (-1)^{2^2 \mu}u
\]
for $u \in V_{\mu}^{(a)}$ with $\mu = 1,0,\frac{1}{2^2}$, preserves the restriction of the algebra product to $V_+^{(a)}$ (i.e. $u^{\sigma(a)} \cdot v^{\sigma(a)} = (u \cdot v)^{\sigma(a)}$ for all $u,v \in V_+^{(a)}$).
\end{description}
\begin{defn}
\label{defn:majorana}
The elements of $A$ are called \emph{Majorana axes} while the automorphisms $\tau(a)$ are called \emph{Majorana involutions}. A real, commutative, non-associative algebra $(V, \cdot )$ equipped with an inner product $( \, , \, )$ is called a \emph{Majorana algebra} if it satisfies axioms M1 and M2 and is generated by a set of Majorana axes.
\end{defn}

The following are basic but important consequences of the Majorana axioms.

\begin{lem}[{\cite[Lemma 1.2]{IPSS10}}]
\label{lem:orthogonality}
Let $V$ be a real vector space satisfying the axioms M1 to M5, and let $a$ be a Majorana axis of $V$. Then the eigenspace decomposition in M4 is orthogonal with respect to $(\, , \, )$ (i.e. $(u,v) = 0$ for all $u \in V_{\mu}^{(a)}$, $v \in V_{\nu}^{(a)}$ with $\mu, \nu \in \{ 1,0,\frac{1}{2^2},\frac{1}{2^5} \}$ and $\mu \neq \nu$).
\end{lem}

\begin{lem}[{\cite[Proposition 2.3.8]{Whybrow18a}}]
\label{lem:conjugate}
Let $V$ be a Majorana algebra generated by a set $A$ of Majorana axes. If $a_0, a_1 \in A$ then the vector $a_0^{\tau(a_1)}$ obeys the axioms M3 - M7 and $\tau(a_0^{\tau(a_1)}) = \tau(a_0)^{\tau(a_1)}$.
\end{lem}

\begin{lem}[{\cite[Lemma 2.2]{IPSS10}}]
\label{lem:conjugate1}
Let $V$ be a Majorana algebra generated by a set $A$ of Majorana axes. If $a_0, a_1 \in A$ then the vector $a_0^{\tau(a_1)}$ is contained in the dihedral algebra $\langle \langle a_0, a_1 \rangle \rangle$.
\end{lem}

In general, if $V$ is a Majorana algebra generated by a set $A$ of Majorana axes and $a_0 \in A$ and $g \in G := \langle \tau(a) \mid a \in A \rangle$ then we usually automatically consider $a_0^g$ to be a Majorana axis and we assume that the set $A$ is closed under the action of $G$. The following lemma gives a similar result for the map $\sigma(a_0)$ in some cases.

\begin{lem}
\label{lem:sigma}
Let $V$ be a Majorana algebra generated by a set $A$ of Majorana axes. If $a_0, a_1 \in A$ such that $\tau(a_1) = 1$ then the vector $a_0^{\sigma(a_1)}$ obeys the axioms M3 - M7 and $\tau(a_0^{\sigma(a_1)}) = \tau(a_0)^{\sigma(a_1)}$.
\end{lem}

\begin{proof}
It is easy to see that $\tau(a_1) = 1$ if and only if $V_{\frac{1}{2^5}}^{(a_1)} = 0$. Thus, from axiom M7, $\sigma(a_1)$ preserves the algebra product on the whole of $V$. In particular, $a_0^{\sigma(a_1)}$ obeys axioms M4 and M5.

Moreover, if $u, v \in V$ then we can write
\[
u = u_1 + u_0  + u_{\frac{1}{2^2}} \textrm{ and } v = v_1 + v_0 + v_{\frac{1}{2^2}}
\]
where $u_{\mu}, v_{\mu} \in V_{\mu}^{(a_1)}$ for $\mu \in \{1, 0, \frac{1}{2^2}\}$. Then, from Lemma \ref{lem:orthogonality},
\[
(u^{\sigma(a_1)}, v^{\sigma(a_1)}) = (u_1, v_1) + (u_0, v_0) + (u_{\frac{1}{2^2}}, v_{\frac{1}{2^2}}) = (u, v).
\]
Thus $\sigma(a_1)$ preserves the inner product and so $a_0^{\sigma(a_1)}$ obeys axiom M3.

We will now show that $a_0^{\sigma(a_1)}$ obeys axioms M6 and M7. In particular, we will show that
\[
\tau(a_0^{\sigma(a_1)}) = \tau(a_0)^{\sigma(a_1)} \textrm{ and } \sigma(a_0^{\sigma(a_1)}) = \sigma(a_0)^{\sigma(a_1)}.
\]
As $\sigma(a_1)$ preserves the algebra product, $v \in V_{\mu}^{(a_0^{\sigma(a_1)})}$ if and only if $v^{\sigma(a_1)} \in V_{\mu}^{(a_0)}$ for $\mu \in \{1, 0, \frac{1}{2^2}, \frac{1}{2^5}\}$. Thus if $v \in V_+^{(a_0^{\sigma(a_1)})}$ then
\[
v^{\sigma(a_1)\tau(a_0)\sigma(a_1)} = (v^{\sigma(a_1)})^{\tau(a_0)\sigma(a_1)} = v^{\sigma(a_1)\sigma(a_1)} = v
\]
and if $v \in V_{\frac{1}{2^5}}^{(a_0^{\sigma(a_1)})}$ then
\[
v^{\sigma(a_1)\tau(a_0)\sigma(a_1)} = (v^{\sigma(a_1)})^{\tau(a_0)\sigma(a_1)} = -v^{\sigma(a_1)\sigma(a_1)} = -v.
\]
Thus $\tau(a_0^{\sigma(a_1)}) = \tau(a_0)^{\sigma(a_1)}$ on $V$. Similarly, we can show that $\sigma(a_0^{\sigma(a_1)}) = \sigma(a_0)^{\sigma(a_1)}$ on $V_+^{(a_0^{\sigma(a_1)})}$. In particular, $\tau(a_0^{\sigma(a_1)})$ and $\sigma(a_0^{\sigma(a_1)})$ are composed of maps which preserve the algebra product on $V$ and $V_+^{(a_0^{\sigma(a_1)})}$ respectively and so $a_0^{\sigma(a_1)}$ obeys axioms M6 and M7.
\end{proof}

The seminal result in Majorana theory was the classification of dihedral Majorana algebras, i.e. Majorana algebras generated by two Majorana axes. This result was inspired by work of Sakuma \cite{Sakuma07} and truly forms the foundations of the theory.

\begin{thm}[{{\cite{IPSS10}}}]
\label{thm:IPSS10}
Let $V$ be a Majorana algebra generated by two Majorana axes $a_0$ and $a_1$. Let $\tau_0 := \tau(a_0)$, $\tau_1 := \tau(a_1)$, $\rho := \tau_0 \tau_1$ and $G := \langle \tau_0, \tau_1 \rangle$. Then
\begin{enumerate}
\item $V$ is isomorphic to a dihedral algebra of type $NX$ for some $1 \le N \le 6$ and $X \in \{A,B,C\}$, the structure of which is given in Table \ref{tab:sakuma};
\item $G \cong \frac{D_{2N}}{Z(D_{2N})}$ for $N$ as above;
\item for $i \in \mathbb{Z}$ and $\epsilon \in \{0,1\}$, the image of $a_{\epsilon}$ under the $i$-th power of $\rho$, which we denote $a_{2i+\epsilon}$, is a Majorana axis and $\tau(a_{2i + \epsilon}) = \rho^{-i}\tau_{\epsilon}\rho^i$.
\end{enumerate}
\end{thm}

\begin{table}
\begin{center}
\vspace{0.35cm}
\noindent
\begin{tabular}{|c|c|c|}
\hline
&&\\
 Type & Basis & Products and angles \\
&&\\
\hline
&&\\
&& $a_0 \cdot a_1=\frac{1}{2^3}(a_0+a_1-a_{\rho}),~a_0 \cdot a_{\rho}=\frac{1}{2^3}(a_0+a_{\rho}-a_1)$ \\
2A & $a_0,a_1,a_{\rho}$ & $a_{\rho} \cdot a_{\rho} = a_{\rho}$ \\
&&$(a_0,a_1)=(a_0,a_{\rho})=(a_{\rho},a_{\rho})=\frac{1}{2^3}$\\
&& \\
2B & $a_0,a_1$ &$a_0 \cdot a_1=0$,~$(a_0,a_1)=0$ \\
&&\\
&  &$a_0 \cdot a_1=\frac{1}{2^5}(2a_0+2a_1+a_{-1})-\frac{3^3 \cdot 5}{2^{11}}u_{\rho}$\\
3A& $a_{-1},a_0,a_1,$ & $a_0 \cdot u_{\rho}=\frac{1}{3^2}(2a_0-a_1-a_{-1})+\frac{5}{2^5}u_{\rho}$~~~~\\
&$u_{\rho}$& $u_{\rho} \cdot u_{\rho}=u_{\rho}$\\
&& $(a_0,a_1)=\frac{13}{2^8}$,~$(a_0,u_{\rho})=\frac{1}{2^2}$,~$(u_{\rho},u_{\rho})=\frac{2^3}{5}$
\\
&&\\
3C & $a_{-1},a_0,a_1$ & $a_0 \cdot a_1=\frac{1}{2^6}(a_0+a_1-a_{-1}),~(a_0,a_1)=\frac{1}{2^6}$ \\
&&\\
&  & ~$a_0 \cdot a_1=\frac{1}{2^6}(3a_0+3a_1+a_2+a_{-1}-3v_{\rho})$\\
4A & $a_{-1},a_0,a_1,$ & $a_0 \cdot v_{\rho}=\frac{1}{2^4}(5a_0-2a_1-a_2-2a_{-1}+3v_{\rho})$\\
&$a_2,v_{\rho}$&~$v_{\rho} \cdot v_{\rho}=v_{\rho}$, ~$a_0 \cdot a_2=0$ \\
& & $(a_0,a_1)=\frac{1}{2^5},~(a_0,a_2)=0,~(a_0,v_{\rho})=\frac{3}{2^3},~(v_{\rho},v_{\rho})=2$\\
&&\\
4B & $a_{-1},a_0,a_1,$ & $a_0 \cdot a_1=\frac{1}{2^6}(a_0+a_1-a_{-1}-a_2+a_{\rho})$
\\
& $a_2,a_{\rho^2}$ & $a_0 \cdot a_2=\frac{1}{2^3}(a_0+a_2-a_{\rho^2})$\\
&& $(a_0,a_1)=\frac{1}{2^6},~(a_0,a_2)=(a_0,a_{\rho^2})=\frac{1}{2^3}$ \\
&&\\
&& $a_0 \cdot a_1=\frac{1}{2^7}(3a_0+3a_1-a_2-a_{-1}-a_{-2})+w_{\rho}$
\\
 5A & $a_{-2},a_{-1},a_0,$ & $a_0 \cdot a_2=\frac{1}{2^7}(3a_0+3a_2-a_1-a_{-1}-a_{-2})-w_{\rho}$
\\
& $a_1,a_2,w_{\rho}$ & $a_0 \cdot w_{\rho}=\frac{7}{2^{12}}(a_{1}+a_{-1}-a_2-a_{-2})+\frac{7}{2^5}w_{\rho}$\\
& & $w_{\rho} \cdot w_{\rho}=\frac{5^2 \cdot 7}{2^{19}}(a_{-2}+a_{-1}+a_0+a_1+a_2)$\\
&&$(a_0,a_1)=\frac{3}{2^7},~(a_0,w_{\rho})=0$, $(w_{\rho},w_{\rho})=\frac{5^3 \cdot 7}{2^{19}}$\\
&& \\
& & $a_0 \cdot a_1=\frac{1}{2^6}(a_0+a_1-a_{-2}-a_{-1}-a_2-a_3+a_{\rho^3})+\frac{3^2 \cdot 5}{2^{11}}u_{\rho^2}$\\
6A& $a_{-2},a_{-1},a_0,$ &$a_0 \cdot a_2=\frac{1}{2^5}(2a_0+2a_2+a_{-2})-\frac{3^3 \cdot 5}{2^{11}}u_{\rho^2}$  \\
&$a_1,a_2,a_3$  &$a_0 \cdot u_{\rho^2}=\frac{1}{3^2}(2a_0-a_2-a_{-2})+\frac{5}{2^5}u_{\rho^2}$  \\
&$a_{\rho^3},u_{\rho^2}$ & $a_0 \cdot a_3=\frac{1}{2^3}(a_0+a_3-a_{\rho^3})$, $a_{\rho^3} \cdot u_{\rho^2}=0$\\
&&$(a_{\rho^3},u_{\rho^2})=0$, $(a_0,a_1)=\frac{5}{2^8}$, $(a_0,a_2)=\frac{13}{2^8}$, $(a_0,a_3)=\frac{1}{2^3}$\\
&&\\
\hline
\end{tabular}
\caption{The dihedral Majorana algebras}
\label{tab:sakuma}
\end{center}
\end{table}
Table \ref{tab:sakuma} does not show all pairwise algebra and inner products of the basis vectors. Those that are missing can be recovered from the action of the group $\langle \tau_0, \tau_1 \rangle$ together with the symmetry between $a_0$ and $a_1$. We also note that the dihedral algebra of type $1A$ is a $1$-dimensional algebra generated by one Majorana axis and so is omitted from Table \ref{tab:sakuma}.

We can use the values in Table \ref{tab:sakuma} to calculate the eigenvectors of these algebras with respect to the axis $a_0$.
\begin{prop}
The eigenspace decompositions with respect to the axis $a_0$ for each dihedral Majorana algebra are given in Table \ref{tab:evecs} (reproduced from \cite[Table 4]{IPSS10}). In each case, the $1$-eigenspace is the $1$-dimensional space spanned by $a_0$ and so is omitted from this table.
\end{prop}
\begin{table}
\begin{center}
\vspace{0.35cm}
\noindent
\begin{tabular}{|>{$}c<{$}|>{$}c<{$}|>{$}c<{$}|>{$}c<{$}|}
\hline
&&&\\
 \mbox{Type} & 0 & \frac{1}{2^2} & \frac{1}{2^5} \\
&&&\\
\hline
&&&\\
2A & a_1 + a_{\rho} - \frac{1}{2^2}a_0 & a_1 - a_{\rho} & \\
&&&\\
2B & a_1 & & \\
&&&\\
3A  & u_{\rho} - \frac{2 \cdot 5}{3^3}a_0 + \frac{2^5}{3^3}(a_1 + a_{-1})
    & \ml{c}{u_{\rho} - \frac{2^3}{3^2 \cdot 5}a_0\\
        - \frac{2^5}{3^2 \cdot 5}(a_1 + a_{-1}) }
    & a_1 - a_{-1} \\
&&&\\
3C  & a_1 + a_{-1} - \frac{1}{2^5}a_0
    &
    & a_1 - a_{-1} \\
&&&\\
4A  & v_{\rho} - \frac{1}{2}a_0 + 2(a_1 + a_{-1}), a_2
    & \ml{c}{v_{\rho} - \frac{1}{3}a_0 \\
        - \frac{2}{3}(a_1 + a_{-1}) - \frac{1}{3}a_2}
    & a_1 - a_{-1} \\
&&&\\
4B  & \ml{c}{a_1 + a_{-1} - \frac{1}{2^5}a_0 - \frac{1}{2^3}(a_{\rho^2}- a_2), \\
        a_2 + a_{\rho^2} - \frac{1}{2^2}a_0}
    & a_2 - a_{\rho^2}
    & a_1 - a_{-1} \\
&&&\\
5A  & \ml{c}{w_{\rho} + \frac{3}{2^9}a_0 - \frac{3 \cdot 5}{2^7}(a_1 + a_{-1}) - \frac{1}{2^7}(a_2 - a_{-2}), \\
        w_{\rho} - \frac{3}{2^9}a_0 + \frac{1}{2^7}(a_1 + a_{-1}) + \frac{3 \cdot 5}{2^7}(a_2 + a_{-2})}
    & \ml{c}{w_{\rho} + \frac{1}{2^7}(a_1 + a_{-1}) \\
        - \frac{1}{2^7}(a_2 + a_{-2})}
    & \ml{c}{a_1 - a_{-1}, \\ a_2 - a_{-2}} \\
&&&\\
6A  & \ml{c}{u_{\rho^2} + \frac{2}{3^2 \cdot 5}a_0 - \frac{2^8}{3^2 \cdot 5}(a_1 - a_{-1}) \\
        - \frac{2^5}{3^2 \cdot 5}(a_2 + a_{-2} + a_3 - a_{\rho^3}), \\
        a_3 + a_{\rho^3} - \frac{1}{2^2}a_0, \\
        u_{\rho^2} - \frac{2 \cdot 5}{3^3}a_0 + \frac{2^5}{3^3}(a_2 + a_{-2})}
    & \ml{c}{u_{\rho^2} - \frac{2^3}{3^2 \cdot 5}a_0  \\
        - \frac{2^5}{3^2 \cdot 5}(a_2 + a_{-2} + a_3)\\
        + \frac{2^5}{3^2 \cdot 5} a_{\rho^3} , \\
        a_3 - a_{\rho^3}}
    & \ml{c}{a_1 - a_{-1}, \\ a_2 - a_{-2}} \\
&&& \\
\hline
\end{tabular}
\caption{The eigenspace decomposition of the dihedral Majorana algebras}
\label{tab:evecs}
\end{center}
\end{table}

The following result follows directly from the values in Table \ref{tab:sakuma}.

\begin{lem}
\label{lem:inclusions}
Let $V$ be a dihedral Majorana algebra (as in Table \ref{tab:sakuma}) that is generated by axes $a_0$ and $a_1$. Then $V$ contains no proper, non-trivial subalgebras, with the exception of the following cases.
\begin{enumerate}[(i)]
\item If $V$ is of type $4A$ or $4B$ then the subalgebras $\langle \langle a_0, a_2 \rangle \rangle$ and $\langle \langle a_1, a_{-1} \rangle \rangle$ are of type $2B$ or $2A$ respectively.
\item If $V$ is of type $6A$ then the subalgebras $\langle \langle a_0, a_3 \rangle \rangle$ and $\langle \langle a_1, a_{-2} \rangle \rangle$ are of type $3A$ and the subalgebras $\langle \langle a_0, a_2 \rangle \rangle$, $\langle \langle a_1, a_{-1} \rangle \rangle$, $\langle \langle a_0, a_{-2} \rangle \rangle$ and $\langle \langle a_1, a_3 \rangle \rangle$ are of type $2A$.
\end{enumerate}
\end{lem}

Informally, this means that we have the following \emph{inclusions} of algebras:
\[
2A \hookrightarrow 4B, \qquad 2B \hookrightarrow 4A, \qquad 2A \hookrightarrow 6A, \qquad \textrm{and} \qquad 3A \hookrightarrow 6A
\]
and that these are the only possible inclusions of non-trivial algebras.

Theorem \ref{thm:IPSS10} shows that if $V$ is a Majorana algebra then there is a finite number of possibilities for the isomorphism types of the dihedral algebras that are contained in $V$. This observation is formalised in the following definition.

\begin{defn}
if $V$ is a Majorana algebra with axes $A$ and $G = \langle \tau(a) \mid A \rangle$ then we define a map $\Psi$ that sends $(a_0, a_1) \in A \times A$ to the type of the dihedral algebra $\langle \langle a_0, a_1 \rangle \rangle$. Then the \emph{shape} of $V$ is the tuple $[\Psi((a_{i_1}, a_{j_1})), \Psi((a_{i_2}, a_{j_2}))  \ldots , \Psi((a_{i_n}, a_{j_n}))]$ where the $(a_{i_k}, a_{j_k})$ are representatives of the orbits of $G$ on $A \times A$.
\end{defn}

\begin{rmk}
\label{rmk:shape}
Suppose that $V$ is a Majorana algebra with Majorana axes $A$ and shape \\ $[\Psi((a_{i_1}, a_{j_1})), \Psi((a_{i_2}, a_{j_2}))  \ldots , \Psi((a_{i_n}, a_{j_n}))]$ where $a_{i_k}, a_{j_k} \in A$ for $1 \leq i \leq n$. Suppose that there exists a permutation $\pi \in \mathrm{Sym}(A)$ that preserves the map $\tau$, i.e. $\tau(a^\pi) = \tau(a)$ for all $a \in A$. Then $V$ must be isomorphic to an algebra with shape $[\Psi((a_{i_1}^\pi, a_{j_1}^\pi)), \Psi((a_{i_2}^\pi, a_{j_2}^\pi))  \ldots , \Psi((a_{i_n}^\pi, a_{j_n}^\pi))]$. In particular, when classifying Majorana algebras, two different shapes give rise to isomorphic algebras if one is the image of the other under such a permutation.
\end{rmk}

Finally, we present the following useful lemma.

\begin{lem}
\label{lem:furtheraxes}
Suppose that $V$ is a Majorana algebra with Majorana axes $A$. Suppose that $a_0, a_1 \in A$ such that $\tau(a_0) = 1$ and such that the subalgebra $U = \langle \langle a_0, a_1 \rangle \rangle$ is of type $2A$. Then the third basis vector $a_{\rho}$ of $U$ obeys the axioms M3 - M7.
\end{lem}

\begin{proof}
From Lemma \ref{lem:sigma}, as $\tau(a_0) = 1$, the map $\sigma(a_0)$ is an automorphism of $V$ and $a_1^{\sigma(a_0)}$ obeys the axioms M3 - M7. From Table \ref{tab:evecs}, we have $a_1 = \frac{1}{2}(\alpha + \beta) + \frac{1}{2^3}a_0$ where
\[
\alpha = a_1 + a_{\rho} - \frac{1}{2^2}a_0 \in V_0^{(a_0)} \textrm{ and } \beta = a_1 - a_{\rho} \in V_{\frac{1}{2^2}}^{(a_0)}.
\]
We then see that
\[
a_1^{\sigma(a_0)} = \frac{1}{2}(\alpha - \beta) + \frac{1}{2^3}a_0 = a_{\rho}
\]
and so $a_{\rho}$ obeys the axioms M3 - M7 as required.
\end{proof}

\subsection{Minimal $3$-generated Majorana algebras}
\label{sec:constructions}

In an analogous way to the situation with groups, we introduce the notion of a minimal $3$-generated Majorana algebra.

\begin{defn}
Let $V$ be a Majorana algebra with Majorana axes $A$. Then $V$ is \emph{minimal $3$-generated} if
\begin{enumerate}[(i)]
\item $V$ is generated by three elements of $A$;
\item if $U$ is a subalgebra of $V$ such that $U = \langle \langle U \cap A \rangle \rangle$ then either $U = V$, or $U$ can be generated by two elements of $A$.
\end{enumerate}
\end{defn}

We will use the list of minimal $3$-generated $6$-transposition groups from Section \ref{sec:groups} in order to partially classify minimal $3$-generated Majorana algebras.

\begin{lem}
\label{lem:min3gen}
Suppose that $V = \langle \langle a_0, a_1, a_2 \rangle \rangle$ is a minimal $3$-generated Majorana algebra such that the group $H = \langle \tau(a_0), \tau(a_1), \tau(a_2) \rangle$ is finite. Then there exists a set $A \subseteq V $ of Majorana axes such that $V = \langle \langle A \rangle \rangle$ and such that $G = \langle \tau(a) \mid a \in A \rangle$ is a minimal $3$-generated $6$-transposition group.
\end{lem}

\begin{proof}
From Lemma \ref{lem:conjugate}, we can take $B := a_0^H \cup a_1^H \cup a_2^H$ to be a set of Majorana axes of $V$. We let
\[
D := \{ \tau(a) \mid a \in B \} = \tau(a_0)^H \cup \tau(a_1)^H \cup \tau(a_2)^H.
\]
We will show that $(H,D)$ is a $3$-generated $6$-transposition group. It is clear that $H$ is generated by three elements of $D$. Moreover, if $t_0, t_1 \in D$, then there exist Majorana axes $b_0, b_1 \in B$ such that $t_i = \tau(b_i)$ for $i \in \{0,1\}$. Then, from Theorem \ref{thm:IPSS10}, we must have $|t_0t_1| = |\tau(b_0) \tau(b_1)| \leq 6$ and so $(H,D)$ is a $6$-transposition group.

If $H$ is a minimal $3$-generated $6$-transposition group then we are done and we set $G = H$. Otherwise, we suppose that $H$ is not minimal $3$-generated. Then Lemma~\ref{l:criterion} implies that there exists $t_0, t_1, t_2 \in D$ such that $H_0 = \langle t_0, t_1, t_2 \rangle$ is a proper non-dihedral subgroup of $H$.

There then exists $b_0,  b_1, b_2 \in B$ such that $\tau(b_i) = t_i$ for $i \in \{0, 1, 2\}$. As $H_0$ is not a dihedral group, $\langle \langle b_0, b_1, b_2 \rangle \rangle$ is not a dihedral algebra and so, as $V$ is minimal $3$-generated, we must have $V = \langle \langle b_0, b_1, b_2 \rangle \rangle$.

Thus if $H_0$ is minimal $3$-generated then we are done and we set $G = H_0$. Otherwise, we continue inductively and construct a family $\mathcal{H} := \{H_i\}_{i = 0}^{\infty}$ of subgroups such that for all $i \geq 0$:
\begin{itemize}
\item $(H_i, H_i \cap D)$ is a $3$-generated $6$-transposition group;
\item $H_{i+1} \leq H_i$;
\item there exists a subset $B_i \leq B$ such that $V = \langle \langle B_i \rangle \rangle$ and $H_i = \langle \tau(a) \mid a \in B_i \rangle$.
\end{itemize}

As $H$ is finite, $\mathcal{H} = \{H_i\}_{i = 0}^{k}$ for some positive integer $k$. Then the group $H_k$ will be a minimal $3$-generated $6$-transposition group and we can take $G = H_k$.
\end{proof}

We note that the converse to this result is not necessarily true. That is to say, it is possible that there exists a Majorana algebra $V$ with axes $A$ such that $G = \langle \tau(a) \mid a \in A \rangle$ is a minimal $3$-generated $6$-transposition group but such that $V$ is not minimal $3$-generated.

This result allows us to use the list of minimal $3$-generated $6$-transposition groups to begin the classification of minimal $3$-generated Majorana algebras. In particular, given a minimal $3$-generated group $G$, we can put strong restrictions on the possibilities for the map $\tau$ from the axes of a Majorana algebra to the involutions of the group. In order to do so, we will require the following results.

\begin{lem}
\label{lem:stabiliser}
Suppose that $V$ is a Majorana algebra with Majorana axes $A$. Let $G = \langle \tau(a) \mid a \in A \rangle$ and write $G_a = \{g \in G \mid a^g = a\}$ where $a \in A$. If $a \in A$ then $\tau(a) \in Z(G_a)$ or, equivalently, $G_a \leq C_G(\tau(a))$.
\end{lem}

\begin{proof}
It is clear that $\tau(a) \in G_a$. If $g \in G_a$ then
\[
[\tau(a), g] = \tau(a)^{-1} \tau(a)^g = \tau(a)\tau(a^g) = \tau(a)^2 = 1.
\]
\end{proof}

\begin{lem}
\label{lem:stabiliser1}
Suppose that $V$ is a Majorana algebra with Majorana axes $A$ and let $G = \langle \tau(a) \mid a \in A \rangle$. If $a \in A$ and $g \in G$ then $G_{a^g} = G_a^g$.
\end{lem}

\begin{proof} For $a$, $g$ and $G_a$ as above,
\begin{align*}
h \in G_{a^g} &\Leftrightarrow a^{gh} = a^g \\
& \Leftrightarrow a^{ghg^{-1}} = a \\
& \Leftrightarrow \exists k \in G_a \textrm{ such that } h = k^g \\
& \Leftrightarrow h \in G_a^g.
\end{align*}
\end{proof}

\begin{lem}
\label{lem:shape}
Suppose that $V$ is a Majorana algebra with Majorana axes $A$. If $a_0, a_1 \in A$ and $D = \langle \tau(a_0), \tau(a_1) \rangle$ then
\begin{enumerate}[(i)]
\item $n = |a_0^D| = |a_1^D|$;
\item if $a_0^D = a_1^D$ then $n \in \{1, 3, 5\}$;
\item if $a_0^D \neq a_1^D$ then $n \in \{1, 2, 3\}$.
\end{enumerate}
Moreover, the subalgebra $\langle \langle a_0, a_1 \rangle \rangle$ is a dihedral algebra of type $NX$ where $N = |a_0^D  \cup a_1^D|$.
\end{lem}

\begin{proof}
This follows directly from the structure of the dihedral algebras given in Table \ref{tab:sakuma}.
\end{proof}

From this result comes the following definition.

\begin{defn}
\label{defn:admissible}
Suppose that $(G, D)$ is a $6$-transposition group and suppose that $V$ is a Majorana algebra with Majorana axes $A$. If $\tau: A \rightarrow D$ is a map then we say that $\tau$ is \emph{admissible} if for all $a_0, a_1 \in A$,
\begin{enumerate}
\item $n := |a_0^H| = |a_1^H|$ where $H = \langle \tau(a_0), \tau(a_1) \rangle$;
\item if $a_0^H= a_1^H$ then $n \in \{1, 3, 5\}$;
\item if $a_0^H \neq a_1^H$ then $n \in \{1, 2, 3\}$.
\end{enumerate}
\end{defn}

Lemma \ref{lem:shape} clearly implies that if $V$ is a Majorana algebra with axes $A$ then the map $\tau$ as defined in axiom M6 is admissible.

\begin{lem}
\label{lem:faithful}
Suppose that $V$ is a Majorana algebra with Majorana axes $A$. Then the group $G = \langle \tau(a) \mid a \in A \rangle$ acts faithfully on $A$.
\end{lem}

\begin{proof}
Suppose that there exists $g \in G$ such that $a^g = a$ for all $a \in A$. As $A$ generates $V$ as an algebra, every element of $A$ can be written as a linear combination of products of the elements of $A$. In particular, these elements must also be fixed by $g$ and so $g = 1$.
\end{proof}

Take $(G,D)$ to be a minimal $3$-generated $6$-transposition group and choose $t_0, t_1, t_2 \in D$ such that $G = \langle t_0, t_1, t_2 \rangle$. We will classify all possible maps $\tau: A \rightarrow G$ such that $V$ is a Majorana algebra with Majorana axes $A$ and such that $V =  \langle \langle a_0, a_1, a_2 \rangle \rangle$ where $a_i \in A$ and $t_i = \tau(a_i)$ for $i \in \{0, 1, 2\}$.

From Lemma \ref{lem:conjugate}, we can take $A = a_0^G \cup a_1^G \cup a_2^G$ and so, using Lemma \ref{lem:faithful}, we can assume that $G$ acts faithfully on $A$ with at most three orbits. Similarly, we need not classify possible maps $\tau$ as above for all triples of involutions that generate $G$. Instead, we consider all triples of conjugacy classes $C_0, C_1, C_2$ of involutions such that $C_0 \cup C_1 \cup C_2$ generates $G$ and take $t_i$ to be a representative of $C_i$ for $i \in \{0, 1, 2\}$.

The action of $G$ on $a_i$ for $i \in \{0, 1, 2\}$ is determined by the orbit $a_i^G$ which is in turn determined by the stabiliser $G_{a_i}$. That is to say, there is a well-defined bijection from the right cosets of $G_{a_i}$ in $G$ to $a_i^G$ given by $G_{a_i}g \mapsto a_i^g$.

From Lemma \ref{lem:stabiliser}, for $i \in \{0, 1, 2\}$, $\langle \tau(a_i) \rangle \leq G_{a_i} \leq C_G(\tau(a_i))$. Moreover, from Lemmas \ref{lem:conjugate} and Lemma \ref{lem:stabiliser1}, if $g \in C_G(\tau(a_i))$ then $\tau(a_i^g) = \tau(a_i)^g = \tau(a_i)$ and $G_{a_i^g} = G_{a_i}^g$. In particular, two stabilisers that are conjugate in $C_G(\tau(a_i))$ lead to the same action and correspond only to taking different representatives of the conjugacy class $C_i$ above.

From this discussion, for a given minimal $3$-generated $6$-transposition group, we can give a method from which it is possible to classify all possible maps $\tau: A \rightarrow G$ where $A$ is a set of Majorana axes of a Majorana algebra $V$ and $\tau$ is the map given in axiom M6.

We note that the orbits $a_0^G$, $a_1^G$ and $a_2^G$ above are not necessarily distinct. In the following method, $k$ denotes the number of orbits of $G$ on $A$ where, in our case, we require that $k \in \{1, 2, 3\}$.

\begin{enumerate}
\item For each $k$ (not necessarily distinct) conjugacy classes of involutions $C_0, \ldots, C_{k-1}$ such that $G$ is generated by $ C_0 \cup \dots \cup C_{k-1}$, choose a representative $t_i$ of $C_i$ for $i \in \{0, \dots , k -1 \}$.
\item Up to conjugacy, find all $k$-tuples $(H_i)_{i = 0}^k$ such that $\langle t_i \rangle \leq H_i \leq C_G(t_i)$ for each $i \in \{0, \ldots, k - 1\}$.
\item Each such tuple $(H_i)_{i = 0}^k$ determines an action of $G$ on $A$ and the length $l_i$ of the orbit of $G$ on the axis $a_i$ will be equal to $[G : H_i]$ for $i \in \{0, \ldots, k - 1\}$.
\item For each tuple $(H_i)_{i = 0}^k$, calculate the permutation homomorphism $\phi_i: G \rightarrow \mathrm{Sym}_{l_i}$ given by the action of right multiplication of $G$ on the set of right cosets $H_i \backslash G$. This then gives the action of $G$ on $a_i$ and the action of $G$ on the whole of $A$ is given by $G^{\phi_0} \times \dots \times G^{\phi_{k - 1}}$. There is then a natural embedding $\phi: G \rightarrow \mathrm{Sym}_l$ where $l = l_0 + \dots + l_{k - 1}$.
\item If $G \not\cong G^\phi $ (i.e. if $G^\phi$ is strictly smaller than $G$) then the permutations giving the action of $G$ on $A$ do not generate the whole of $G$ and so the tuple $(H_i)_{i = 0}^k$ does not give a valid action of $G$ on a set of Majorana axes and we can discard this possibility.
\item Otherwise, we now define a set of Majorana axes $A := a_0^G \cup \dots \cup a_{k-1}^G$ where, for each $i \in \{0, \dots, k - 1\}$, we define that
    \begin{enumerate}[i)]
    \item $\tau(a_i) = t_i$;
    \item $a_i^g = a_{i^{g^\phi}}$;
    \item $\tau(a_i^g) = \tau(a_i)^g$.
    \end{enumerate}
\item Finally, we check that the map $\tau$ is admissible (as in Definition \ref{defn:admissible}). If this is the case, then we can use Lemma \ref{lem:shape} to determine the possible shapes of a Majorana algebra generated by $A$.
\end{enumerate}

\begin{thm}
\label{thm:majoranamain}
Table \ref{tab:majoranamain} gives the shapes and dimensions of all non-zero Majorana algebras $V$ with Majorana axes $A$ such that
\begin{enumerate}[(i)]
\item $V$ is minimal $3$-generated and $V = \langle \langle a_0, a_1, a_2 \rangle \rangle$ for some $a_0, a_1, a_2 \in A$;
\item the group $G = \langle \tau(a_0), \tau(a_1), \tau(a_2) \rangle$ is isomorphic to one of the groups listed in Theorem \ref{t:Min3GenGroups} except possibly $(5^2:3):2$ or $p^2:2$ for $p \in \{3, 5\}$.
\end{enumerate}
\end{thm}

The shapes of each of these Majorana algebras is defined by the information in Tables \ref{tab:shapes3}, \ref{tab:shapes2} and \ref{tab:shapes1}. In the case where the orbits $a_0^G$, $a_1^G$ and $a_2^G$ are distinct  (Table \ref{tab:shapes3}), the columns are labelled as below
\begin{description}
    \item[Group] - the isomorphism type of the group $G = \langle \tau(a_0), \tau(a_1), \tau(a_2) \rangle$;
    \item[Axes] - the lengths of the orbits of $G$ on $A$ written as $|a_0^G| + |a_1^G| + |a_2^G|$;
    \item[No. of shapes] - the number of possible shapes of a Majorana algebra with the given group and axes;
    \item[Orbitals] - representatives of orbits of $G$ on $A \times A$;
    \item[Shape] - for each representative in the column ``orbitals'', a string representative the type of the corresponding dihedral algebra.
\end{description}

If there are just two or one orbits of $G$ on $A$ then the column ``axes'' of Table \ref{tab:shapes2} or \ref{tab:shapes1} gives the value $|a_0^G| + |a_1^G|$ or $|a_0^G|$ respectively. In these cases, we also include an additional column, labelled ``relations'', that allow us to recover the value of $a_1$ and $a_2$ where necessary.

Note that we record values in the columns ``orbitals'' and ``shape'' only where a choice must to be made on the type of the corresponding algebra. We denote by $NX^k$ a tuple of $k$ subalgebras that are of type $NX$ for $X \in \{A, B, C\}$ possibly different for each subalgebra. All other types of dihedral subalgebras are determined by Lemmas \ref{lem:shape} and \ref{lem:inclusions}.

As explained in Remark \ref{rmk:shape}, a number of different shapes may lead to isomorphic algebras. For each such set of shapes, we include just one representative shape and its corresponding algebra in Tables \ref{tab:majoranamain}, \ref{tab:incomplete}, \ref{tab:nonmin3gencomplete} and \ref{tab:nonmin3genincomplete}.

For four of the possible shapes, we have been unable to classify whether the corresponding Majorana algebras are minimal $3$-generated. These cases are given in Table \ref{tab:incomplete}. For many of the possible remaining cases, it turns out that no minimal $3$-generated Majorana algebras of this shape exist. The cases which admit non-minimal $3$-generated Majorana algebras are discussed in Section \ref{sec:nonminimal}.

\begin{proof}
We have implemented the method above as an algorithm in GAP which takes as its input a group and an integer $k \in \{1, 2, 3\}$ and returns the images in $G$ of all possible maps $\tau: A \rightarrow G$ where $A$ is a set of Majorana axes and $G$ acts on $A$ with $k$ orbits. With this information, we can then use the algorithm described in \cite{PW18a} and implemented in \cite{PW18b} to construct the possible algebras satisfying the conditions in the statement of the theorem. It then suffices to check that the algebra in question is indeed minimal $3$-generated.

In all but two cases, this method is sufficient to construct the Majorana algebras for the shape in question. We deal with the two remaining cases separately. The first is the case where $G \cong 2^2$ and where the action of $G$ on $A$ consists of three orbits, each of length two, and leads to two possible shapes. When this shape is chosen to be $(4A)$, the algebra in question exhibits very unusual behaviour; there is an infinite family of Majorana algebras satisfying this description, each with dimension $12$. The full construction of this family is given in \cite{Whybrow18b}.

We finally consider the case where $G \cong 1$ and where $G$ acts on $A$ with three orbits, each of length $1$. We will show that there are two possible algebras of this form with shape $(2A, 2A, 2A)$. We first label the elements of $A$ as $a_0$, $a_1$ and $a_2$ and let $V = \langle \langle A \rangle \rangle$. As  the algebra $V$ has shape $(2A, 2A, 2A)$, it must also contain the three vectors $a_{\rho(0,1)}$, $a_{\rho(0,2)}$ and $a_{\rho(1,2)}$ where $a_{\rho(i,j)}$ denotes the third basis vector from the $2A$ dihedral algebra $\langle \langle a_i, a_j \rangle \rangle$ where $i, j \in \{0, 1, 2\}$ and $i \neq j$.

From Lemma \ref{lem:furtheraxes}, for $i, j \in \{0, 1, 2\}$ such that $i \neq j$, the vector $a_{\rho(i,j)}$ obeys the axioms M3 - M7 and $\tau(a_{\rho(i,j)}) = 1$. In particular, the three algebras $\langle \langle a_i, a_{\rho(j,k)} \rangle \rangle$ for $\{i, j, k\} = \{1, 2, 3\}$ are each isomorphic to a dihedral subalgebra of the Griess algebra of type $2A$ or $2B$. We thus have a futher choice on the types of these algebras.

We use the algorithm implemented in \cite{PW18b} to show that if we choose the algebras $\langle \langle a_i, a_{\rho(j,k)} \rangle \rangle$ to all be either of type $2B$ or of type $2A$ then $V$ is a Majorana algebra of dimension $6$ or $9$ respectively.
\end{proof}

\begin{table}
\renewcommand\arraystretch{0.75}
\begin{center}
\vspace{0.35cm}
\noindent
\begin{tabular}{|>{$}c<{$}  >{$}c<{$} >{$}c<{$} >{$}c<{$} >{$}c<{$}  |} \hline
&&&& \\
G & \textrm{Axes} & \textrm{No. of shapes} & \textrm{Orbitals} & \textrm{Shape} \\
&&&& \\ \hline
&&&& \\
1 & 1 + 1 + 1 & 8 & (a_0, a_1), (a_0, a_2), (a_1, a_2) & 2X^3 \\
&&&& \\
2^2 & 1 + 2 + 2 & 8 & (a_0, a_1), (a_0, a_2), (a_1, a_2) & 2X^2, 4X \\
&&&& \\
2^2 & 2 + 2 + 2 & 8 & (a_0, a_1), (a_0, a_1^{\tau_2}), (a_0, a_2),  & 2X^2, 4X \\
&&&& \\
2^2 & 2 + 2 + 2 & 2 & (a_0, a_1) & 4X \\
&&&& \\
S_3 & 1 + 3 + 3 & 4 & (a_0, a_1), (a_0, a_2) & 2X^2 \\
&&&& \\
S_3 & 3 + 3 + 3 & 1 & & \\
&&&& \\
2^3 &  4 + 4 + 4  & 64 & \renewcommand\arraystretch{1.2}\ml{c}{(a_0, a_1), (a_0, a_2), (a_1, a_2), \\ (a_0, a_0^{\tau_1\tau_2}), (a_1, a_1^{\tau_0\tau_2}), (a_2, a_2^{\tau_0\tau_2})} & 4X^3, 2X^3 \\
&&&& \\
2^3 &  4 + 4 + 2  &  16 & \renewcommand\arraystretch{1.2}\ml{c}{(a_0, a_1), (a_0, a_2), \\ (a_0, a_0^{\tau_1\tau_2}), (a_1, a_1^{\tau_0\tau_2})} & 4X^2, 2X^2 \\
&&&& \\
2^3 &  4 + 2 + 2 & 32 & \renewcommand\arraystretch{1.2}\ml{c}{(a_0, a_1), (a_0, a_2), (a_1, a_2),\\ (a_0, a_0^{\tau_1\tau_2}), (a_1, a_2^{\tau_0}) }& 4X^2, 2X^4 \\
&&&& \\
2^3 &  4 + 2 + 2  &  4 & (a_0, a_1), (a_0, a_0^{\tau_1\tau_2})  & 4X, 2X \\
&&&& \\
D_{12} &   6 + 6 + 6 & 4 & (a_0, a_1), (a_2, a_2^{\tau_0\tau_1}) & 4X, 2X \\
&&&& \\
D_{12} &   6 + 6 + 2 & 2 & (a_0, a_2) &     4X \\
&&&& \\
3^2 : 2 &   9 + 9 + 9 & 1  & & \\
&&&& \\
S_4 &  12 + 12 + 12 & 2 & (a_0, a_1) & 4X \\
&&&& \\
S_4 &  6 + 6 + 6 & 64  & \renewcommand\arraystretch{1.2}\ml{c}{(a_0, a_1), (a_0, a_0^{(\tau_0\tau_1\tau_2)^2}), \\ (a_0, a_2^{(\tau_2\tau_1)^2}),  (a_1, a_1^{(\tau_1\tau_0\tau_2)^2}), \\ (a_1, a_2^{(\tau_2\tau_0)^2}), (a_2, a_2^{(\tau_0\tau_1)})} & 2X^6 \\
&&&& \\
S_4 &  6 + 6 + 3 & 16  & \renewcommand\arraystretch{1.2}\ml{c}{(a_0, a_2), (a_1, a_2), \\(a_0, a_2^{\tau_0\tau_1}), (a_0, a_1^{\tau_2\tau_0\tau_2\tau_1})} & 4X, 2X^3 \\
&&&& \\
(3^3:3):2 &   27 + 27 + 27  & 1  & & \\
&&&& \\
(3^3:3):2 &   27 + 27 + 9 & 1  & & \\
&&&& \\
(3^3:3):2 &   27 + 9 + 9 & 1  & & \\
&&&& \\ \hline
\end{tabular}
\caption{The number of possible shapes of algebras with $3$ orbits of axes}
\label{tab:shapes3}
\end{center}
\end{table}

\begin{table}
\renewcommand\arraystretch{0.75}
\begin{center}
\vspace{0.35cm}
\noindent
\begin{tabular}{|>{$}c<{$}  >{$}c<{$} >{$}l<{$} >{$}c<{$} >{$}c<{$} >{$}c<{$}  |} \hline
&&&&& \\
G & \textrm{Axes} & \textrm{Relations} & \textrm{No. of shapes} & \textrm{Orbitals} & \textrm{Shape} \\
&&&&& \\ \hline
&&&&& \\
1 &  1 + 1 & & 2 & (a_0, a_1) & 2X \\
&&&&& \\
2^2 &  2 + 2 & a_2 = a_1^{\tau_0} & 2 & (a_0, a_1) & 4X \\
&&&&& \\
S_3 &  1 + 3 & a_2 = a_1^{\tau_2\tau_1} & 4 & (a_0, a_1), (a_1, ,a_2) & 2X, 3X \\
&&&&& \\
S_3 &  3 + 3 & a_2 = a_1^{\tau_0\tau_1} & 1 & & \\
&&&&& \\
D_{10} &  1 + 5 & a_2 = a_1^{(\tau_2\tau_1)^2} & 2 & (a_0, a_1) & 2X \\
&&&&& \\
D_{12} &  6 + 6 & a_2 = a_0^{\tau_2\tau_0} & 4 & (a_0, a_1), (a_0, a_1^{\tau_2\tau_0}) & 4X, 2X \\
&&&&& \\
D_{12} &  6 + 2 & a_2 = a_0^{\tau_2\tau_0} & 2 & (a_1, a_1^{\tau_0}) & 2X \\
&&&&& \\
3^2 : 2 &  9 + 9 & a_2 = a_0^{\tau_2\tau_0} & 1 & & \\
&&&&& \\
S_4 &  12 + 12 & a_2 = a_0^{(\tau_0\tau_2\tau_1)^2} & 2 & (a_0, a_2) & 4X \\
&&&&& \\
S_4 &  6 + 6 & a_2 = a_0^{\tau_1\tau_0} & 8 & \renewcommand\arraystretch{1.2}\ml{c}{(a_0, a_0^{(\tau_0\tau_2)^{\tau_1}}), (a_0, a_1), \\ (a_1, a_1^{(\tau_0\tau_2)^{\tau_1\tau_0}})} & 2X^3 \\
&&&&& \\
S_4 &  6 + 3 & a_2 = a_0^{\tau_2\tau_0} &  8 & (a_0, a_1), (a_0, a_2), (a_0, a_1^{\tau_2\tau_0}) & 4X, 3X, 2X \\
&&&&& \\
(3^2 : 3): 2 &  27 + 27 & a_2 = a_0^{\tau_2\tau_0} & 1 & & \\
&&&&& \\
(3^2 : 3): 2 &  27 + 9 & a_2 = a_0^{\tau_2\tau_0} & 1 & & \\
&&&&& \\ \hline
\end{tabular}
\caption{The number of possible shapes of algebras with $2$ orbits of axes}
\label{tab:shapes2}
\end{center}
\end{table}

\begin{table}
\renewcommand\arraystretch{0.75}
\begin{center}
\vspace{0.35cm}
\noindent
\begin{tabular}{|>{$}c<{$}  >{$}c<{$} >{$}l<{$} >{$}c<{$} >{$}c<{$} >{$}c<{$}  |} \hline
&&&&& \\
G & \textrm{Axes} & \textrm{Relations} & \textrm{No. of shapes} & \textrm{Orbitals} & \textrm{Shape} \\
&&&&& \\ \hline
&&&&& \\
1 &  1 & & 1 & & \\
&&&&& \\
S_3 &  3 & a_1 = a_0^{\tau_2}, a_2 = a_0^{\tau_1} & 2 & (a_0, a_1) & 3X \\
&&&&& \\
D_{10} &  5 & a_1 = a_0^{\tau_1\tau_0}, a_2 = a_0^{\tau_1} & 1 &  &  \\
&&&&& \\
3^2 :2 &  9 & a_1 = a_0^{\tau_1\tau_0}, a_2 = a_0^{\tau_2\tau_0}  & 16 & \renewcommand\arraystretch{1.2}\ml{c}{ (a_0, a_1), (a_0, a_2) \\ (a_0, a_2^{\tau_0\tau_1}), (a_0, a_2^{\tau_1\tau_0})} & 3X^4 \\
&&&&& \\
S_4 &  12  & a_1 = a_0^{(\tau_0\tau_1\tau_2)^2}, a_2 = a_0^{\tau_2\tau_0} & 2 & (a_0, a_1) & 4X \\
&&&&& \\
S_4 &   6  & a_1 = a_0^{(\tau_0\tau_1)^{\tau_2}}, a_2 = a_0^{\tau_2\tau_0} & 4 & (a_0, a_1), (a_0, a_2) & 2X, 3X \\
&&&&& \\
5^2 : 2 &  25 & a_1 = a_0^{(\tau_1\tau_0)^2}, a_2 = a_0^{(\tau_2\tau_0)^2} & 1 &  &  \\
&&&&& \\
(3^2 : 3): 2 & 27 & a_1 = a_0^{(\tau_0\tau_1)^{\tau_2}}, a_2 = a_0^{(\tau_0\tau_2)^{\tau_1}} & 1 & &  \\
&&&&& \\
A_5 &  15 & a_1 = a_0^{(\tau_1\tau_0)^{\tau_2}\tau_1\tau_2}, a_2 = a_0^{\tau_2\tau_0} & 4 & (a_0, a_1), (a_0, a_2) & 2X, 3X \\
&&&&& \\
(5^2 : 3) : 2 &  15  & a_1 = a_0^{\tau_1\tau_0}, a_2 = a_0^{\tau_2\tau_0} & 2 & (a_0, a_1) & 3X \\
&&&&& \\ \hline
\end{tabular}
\caption{The number of possible shapes of algebras with $1$ orbit of axes}
\label{tab:shapes1}
\end{center}
\end{table}

\begin{table}
\renewcommand\arraystretch{1.1}
\begin{center}
\vspace{0.35cm}
\noindent
\begin{tabular}{|>{$}c<{$}   >{$}c<{$} >{$}c<{$} >{$}c<{$}  |} \hline
    &&& \\
G & \textrm{Axes}  & \mathrm{Shape} & \textrm{dim } V  \\
&&& \\
 \hline
\rule{0pt}{4ex} 1 & 1 + 1 + 1 & (2B, 2B, 2B) & 3 \\
 & & (2B, 2B, 2A) & 4 \\
 & & (2B, 2A, 2A) & 6 \\
 & & (2A, 2A, 2A) & 6, 9 \\ 
\rule{0pt}{2ex} & 1 + 1 & (2B) & 2 \\
 & & (2A) & 3 \\
\rule{0pt}{2ex} &  1 &  & 1 \\
\rule{0pt}{4ex} 2^2 & 1 + 2 + 2 & (2A, 2A, 4B) & 5 \\
\rule{0pt}{2ex} & 2 + 2 + 2 & (4A) & 12 \\ 
 & & (4B) & 7 \\
\rule{0pt}{2ex} &  2 + 2 & (4A) & 5 \\
 & & (4B) &  5 \\
\rule{0pt}{4ex} S_3 &  1 + 3 & (2B, 3C) & 4 \\
& & (2B, 3A) & 5 \\
& & (2A, 3A) & 8 \\ 
\rule{0pt}{4ex} &  3 + 3 & & 8 \\
\rule{0pt}{2ex} &  3 & (3C) &  3 \\
& &   (3A) &  4 \\
\rule{0pt}{4ex} D_{10} &  1 + 5 & (2B) & 7 \\
\rule{0pt}{2ex} &  5 &  &  6 \\
\rule{0pt}{4ex} D_{12} &  6 + 2 & (4B, 2A) & 10 \\
\rule{0pt}{4ex} 3^2 :2 &  9 &  (3C, 3C, 3C, 3C) &  9 \\
& &  (3C, 3C, 3C, 3A) &  12 \\
\rule{0pt}{4ex} S_4 &  6 + 3 &  (4B, 3C, 2A) & 9  \\
& & (4B, 3A, 2A) & 13  \\
\rule{0pt}{2ex} &  12 &  (4B) &  17 \\
\rule{0pt}{2ex} &   6  &  (2B, 3C) &  6 \\
& &   (2A, 3C) &  9 \\
& &   (2B, 3A) &  13 \\
& &   (2A, 3A) & 13 \\
\rule{0pt}{4ex} A_5 &  15 &  (2A, 3C) &  20 \\
& &   (2A, 3A) &  26 \\
\rule{0pt}{4ex} (5^2 : 3) : 2 &  15 &  (3C) &  18 \\
&&& \\ \hline
\end{tabular}
\caption{Minimal $3$-generated Majorana algebras}
\label{tab:majoranamain}
\end{center}
\end{table}

\begin{table}
\renewcommand\arraystretch{1.1}
\begin{center}
\vspace{0.35cm}
\noindent
\begin{tabular}{|>{$}c<{$}  >{$}c<{$} >{$}c<{$}  |} \hline
&& \\
G & \textrm{Axes}  & \mathrm{Shape}   \\
&& \\ \hline
\rule{0pt}{4ex} 3^2 :2 &  9 &  (3C, 3C, 3A, 3A) \\
& & (3A, 3A, 3A, 3A) \\
\rule{0pt}{4ex} 5^2 : 2 &  25 &  \\
\rule{0pt}{4ex} (5^2 : 3) : 2 &  15 &  (3A) \\
&& \\\hline
\end{tabular}
\caption{Unresolved cases}
\label{tab:incomplete}
\end{center}
\end{table}

\subsection{The non-minimal examples}
\label{sec:nonminimal}

Recall from the discussion following Lemma \ref{lem:min3gen} it is possible that a minimal $3$-generated $6$-transposition group admits a Majorana algebra that is not itself minimal $3$-generated. In Tables \ref{tab:nonmin3gencomplete} and \ref{tab:nonmin3genincomplete}, we list all cases where our method in Section \ref{sec:constructions} gives rise to non-zero Majorana algebras that are not minimal $3$-generated. In each case, we give details of the original algebra and also of the minimal $3$-generated subalgebras it contains.

Table \ref{tab:nonmin3genincomplete}, gives the examples where the original algebra has not been fully constructed. In these cases, we can still conclude that these algebras are not minimal $3$-generated but it is possible that the list of minimal $3$-generated subalgebras is not exhaustive.

We note that in all but one case the minimal $3$-generated subalgebras have been found from other groups that we have considered. The exception is the case where $G \cong 3^2 : 2$ acts with two orbits each of length nine on the Majorana axes of a Majorana algebra $V$. In this case $V$ is not minimal $3$-generated and contains an algebra of dimension $21$. This must arise from a group of the same isomorphism type acting with one orbit of length nine on the Majorana axes of a Majorana algebra. This algebra must have shape involving only $3A$, rather than $3C$, algebras. We have not been able to directly classify algebra of this form. This strongly suggests that there are at least two algebras with this shape.

\section*{Acknowledgements} We would like to express our deepest thanks to Prof. S.~Shpectorov for his support throughout this project. His comments and advice have been extremely helpful to us in the writing of this paper. We would also like to thank an anonymous referee for their detailed and constructive comments.

\pagebreak

%
%

\renewcommand\arraystretch{0.9}
\begin{longtable}[c]{|>{$}c<{$}>{$}c<{$}>{$}c<{$}>{$}c<{$}|>{$}c<{$}>{$}c<{$}>{$}c<{$}>{$}c<{$}|} \hline
\multicolumn{4}{|c|}{ Original algebra } & \multicolumn{4}{c|}{ Minimal $3$-generated subalgebras} \\ \hline
&&&&&&& \\
G & \mathrm{Axes} & \mathrm{Shape} & \textrm{dim } V & G & \mathrm{Axes} & \mathrm{Shape} & \textrm{dim } V \\
&&&&&&& \\ \hline
\endfirsthead
\endhead
\endfoot
\endlastfoot
&&&&&&& \\
2^2 & 1 + 2 + 2 &  ( 2B, 2B, 4A) & 6 & 1 & 1 + 1 + 1 & (2B, 2B, 2B) & 3 \\
&&&&&&& \\
2^2 & 1 + 2 + 2 &  ( 2B, 2B, 4B) & 6 & 1 & 1 + 1 + 1 & (2B, 2B, 2A) & 4 \\
&&&&&&& \\
2^2 & 1 + 2 + 2 &  ( 2B, 2A, 4A) & 10 & 1 & 1 + 1 + 1 & (2B, 2B, 2B) & 3 \\
& & & & 1 & 1 + 1 + 1 & (2B, 2A, 2A) & 6 \\
&&&&&&& \\
2^2 & 1 + 2 + 2 &  ( 2B, 2A, 4B) & 8 & 1 & 1 + 1 + 1 & (2B, 2B, 2A) & 4 \\
& & & & 1 & 1 + 1 + 1 & (2B, 2A, 2A) & 6 \\
&&&&&&& \\
2^2 & 1 + 2 + 2 &  ( 2A, 2A, 4A) & 14 & 1 & 1 + 1 + 1 & (2B, 2A, 2A) & 6 \\
&&&&&&& \\
2^2 & 2 + 2 + 2 &  ( 2B, 2B, 4A ) & 9 & 1 & 1 + 1 + 1 & (2B, 2B, 2B) & 3 \\
&&&&&&& \\
2^2 & 2 + 2 + 2 &  ( 2B, 2A, 4B ) & 8 & 1 & 1 + 1 + 1 & (2B, 2A, 2A) & 6 \\
&&&&&&& \\
2^2 & 2 + 2 + 2 &  ( 2A, 2A, 4B ) & 11 & 1 & 1 + 1 + 1 & (2A, 2A, 2A) & 9 \\
&&&&&&& \\
S_3 & 1 + 3 + 3 &  ( 2B, 2B ) & 9 & S_3 & 1 + 3 & (2B, 3A) & 5 \\
& & & & 1 & 1 + 1 + 1 & (2B, 2B, 2A) & 4 \\
&&&&&&& \\
S_3 & 1 + 3 + 3 &  ( 2B, 2A ) & 13 & S_3 & 1 + 3 & (2B, 3A) & 5 \\
& & & & 1 & 1 + 1 + 1 & (2B, 2A, 2A) & 6 \\
&&&&&&& \\
S_3 & 1 + 3 + 3 &  ( 2A, 2A ) & 11 & 1 & 1 + 1 + 1 & (2B, 2A, 2A) & 6 \\
&&&&&&& \\
S_3 & 3 + 3 + 3 &  & 13 & 1 & 1 + 1 + 1 & (2B, 2A, 2A) & 6 \\
&&&&&&& \\
2^3  & 4 + 4 + 4 &  ( 4B, 4B, 4B, 2B, 2B, 2B) & 15 & 1 & 1 + 1 + 1 & (2B, 2A, 2A) & 6 \\
&&&&&&& \\
2^3  & 4 + 4 + 2 &  ( 4A, 4B, 2A, 2A ) & 13 & 1 & 1 + 1 + 1 & (2B, 2A, 2A) & 6 \\
&&&&&&& \\
2^3  & 4 + 2 + 2 &  ( 4A, 4A, 2B, 2B, 2B) & 13 & 1 & 1 + 1 + 1 & (2B, 2B, 2B) & 3 \\
&&&&&&& \\
2^3  & 4 + 2 + 2 &  ( 4A, 4B, 2A, 2B, 2B) & 12 & 1 & 1 + 1 + 1 & (2B, 2A, 2A) & 6 \\
& & & & 1 & 1 + 1 + 1 & (2B, 2B, 2B) & 3 \\
& & & & 1 & 1 + 1 + 1 & (2B, 2B, 2A) & 4 \\
&&&&&&& \\
2^3  & 4 + 2 + 2 &  ( 4A, 4B, 2A, 2B, 2A) & 15 & 1 & 1 + 1 + 1 & (2B, 2A, 2A) & 6 \\
& & & & 1 & 1 + 1 + 1 & (2B, 2B, 2A) & 4 \\
&&&&&&& \\
2^3  & 4 + 2 + 2 &  ( 4B, 4B, 2B, 2B, 2B) & 10 & 1 & 1 + 1 + 1 & (2B, 2A, 2A) & 6 \\
& & & & 1 & 1 + 1 + 1 & (2B, 2B, 2A) & 4 \\
&&&&&&& \\
3^2 : 2 & 9 + 9 & & 31 & 3^2:2 & 9 & (3A, 3A, 3A, 3A) & 21 \\
&&&&&&& \\
3^2:2 & 9 + 9 + 9 &  & 42 & 3^2:2 & 9 & (3A, 3A, 3A, 3A) & 21 \\
& & & & 1 & 1 + 1 + 1 & (2B, 2A, 2A) & 6 \\
&&&&&&& \\
S_4 & 6 + 6 & (2B, 2A, 2A) & 20 & S_4  &  6 & (2A, 3A) & 13  \\
& &  &  &  1 &  1 + 1 + 1 & (2B, 2A, 2A) &  6  \\
&&&&&&& \\
S_4 & 6 + 3 & (4A, 3C, 2B) & 12 & S_4 & 6 & (2B, 3C) & 6 \\
& & & & 1 & 1 + 1 + 1 & (2B, 2B, 2B) & 3 \\
&&&&&&& \\
S_4 & 6 + 3 & (4A, 3A, 2B) & 25 & S_4 & 6 & (2B, 3A) & 13 \\
& & & & 1 & 1 + 1 + 1 & (2B, 2B, 2B) & 3 \\
&&&&&&& \\
S_4 & 6 + 3 & (4A, 3A, 2A) & 23 & S_4 & 6 & (2A, 3A) & 13 \\
& & & & 1 & 1 + 1 + 1 & (2B, 2A, 2A) & 6 \\
&&&&&&& \\
S_4 & 6 + 3 & (4B, 3C, 2B) & 12 & S_4 & 6 & (2A, 3C) & 9 \\
& & & & 1 & 1 + 1 + 1 & (2B, 2B, 2A) & 4 \\
&&&&&&& \\
S_4 & 6 + 3 & (4B, 3A, 2B) & 16 & S_4 & 6 & (2B, 3C) & 13 \\
& & & & 1 & 1 + 1 + 1 & (2B, 2B, 2A) & 4 \\
&&&&&&& \\
S_4 & 6 + 6 + 3 &  ( 4B, 2B, 2B, 2B) & 23 & S_4 & 6 & (2B, 3A) & 13 \\
& & & & 1 & 1 + 1 + 1 & (2B, 2A, 2A) & 6 \\
& & & & 1 & 1 + 1 + 1 & (2B, 2B, 2A) & 4 \\
& & & & 1 & 1 + 1 + 1 & (2B, 2B, 2B) & 3 \\
&&&&&&& \\
S_4 & 6 + 6 + 3 &  ( 4B, 2B, 2A, 2A) & 20 & S_4 & 6 & (2B, 3A) & 13 \\
& & & & 1 & 1 + 1 + 1 & (2B, 2A, 2A) & 6 \\
& & & & 1 & 1 + 1 + 1 & (2B, 2A, 2A) & 6 \\
&&&&&&& \\
S_4 & 6 + 6 + 6 &  ( 2A, 2B, 2B, 2A, 2B, 2A) & 28 & S_4 & 6 & (2B, 3A) & 13 \\
& & & & 1 & 1 + 1 + 1 & (2B, 2A, 2A) & 6 \\
& & & & 1 & 1 + 1 + 1 & (2B, 2B, 2A) & 4 \\
&&&&&&& \\
A_5 &   15  &  (2B, 3C) &  21 & 1 & 1 + 1 + 1 & (2B, 2B, 2B) & 3 \\
&&&&&&& \\
A_5 &   15  &  (2B, 3A) &  46 & 1 & 1 + 1 + 1 & (2B, 2B, 2B) & 3 \\
&&&&&&& \\ \hline
\caption{The complete non-minimal examples }
\label{tab:nonmin3gencomplete}
\end{longtable}

\pagebreak

%
%

\renewcommand\arraystretch{0.9}
\begin{longtable}[c]{|>{$}c<{$}>{$}c<{$}>{$}c<{$}|>{$}c<{$}>{$}c<{$}>{$}c<{$}>{$}c<{$}|} \hline
\multicolumn{3}{|c|}{ Original algebra } & \multicolumn{4}{c|}{ Minimal $3$-generated subalgebras} \\ \hline
&&&&&& \\
G & \mathrm{Axes} & \mathrm{Shape} & G & \mathrm{Axes} & \mathrm{Shape} & \textrm{dim } V \\
&&&&&& \\ \hline
\endfirsthead
\endhead
\endfoot
\endlastfoot
&&&&&& \\
2^2 & 2 + 2 + 2 &  ( 2A, 2A, 4A ) & 1 & 1 + 1 + 1 & (2B, 2A, 2A) & 6 \\ 
&&&&&& \\
2^3  & 4 + 4 + 4 &  ( 4A, 4A, 4A, 2B, 2B, 2B)  & 1 & 1 + 1 + 1 & (2B, 2B, 2B) & 3 \\
&&&&&& \\
2^3  & 4 + 4 + 4 &  ( 4A, 4A, 4A, 2B, 2B, 2A)  & 1 & 1 + 1 + 1 & (2B, 2B, 2B) & 3 \\
& & &  1 & 1 + 1 + 1 & (2B, 2B, 2A) & 4 \\
&&&&&& \\
2^3  & 4 + 4 + 4 &  ( 4A, 4A, 4A, 2A, 2B, 2A)  & 1 & 1 + 1 + 1 & (2B, 2B, 2B) & 3 \\
& & &  1 & 1 + 1 + 1 & (2B, 2B, 2A) & 4 \\
&&&&&& \\
2^3  & 4 + 4 + 4 &  ( 4B, 4A, 4A, 2A, 2B, 2A)  & 1 & 1 + 1 + 1 & (2B, 2B, 2B) & 3 \\
& & &  1 & 1 + 1 + 1 & (2B, 2B, 2A) & 6 \\
&&&&&& \\
2^3  & 4 + 4 + 4 &  ( 4A, 4A, 4A, 2A, 2A, 2A)  & 1 & 1 + 1 + 1 & (2B, 2B, 2A) & 4 \\
&&&&&& \\
2^3  & 4 + 4 + 4 &  ( 4B, 4A, 4A, 2A, 2A, 2A)  & 1 & 1 + 1 + 1 & (2B, 2B, 2A) & 4 \\
& & &  1 & 1 + 1 + 1 & (2B, 2A, 2A) & 6 \\
&&&&&& \\
2^3  & 4 + 4 + 4 &  ( 4B, 4A, 4B, 2A, 2A, 2B)  & 1 & 1 + 1 + 1 & (2B, 2A, 2A) & 6 \\
&&&&&& \\
2^3  & 4 + 4 + 4 &  ( 4A, 4B, 4A, 2A, 2A, 2B)  & 1 & 1 + 1 + 1 & (2B, 2A, 2A) & 6 \\
& & &  1 & 1 + 1 + 1 & (2B, 2B, 2B & 3 \\
&&&&&& \\
2^3  & 4 + 2 + 2 &  ( 4A, 4A, 2B, 2B, 2A)  & 1 & 1 + 1 + 1 & (2B, 2B, 2B) & 3 \\
& & &  1 & 1 + 1 + 1 & (2B, 2B, 2A) & 4 \\
&&&&&& \\
2^3  & 4 + 2 + 2 &  ( 4A, 4A, 2B, 2A, 2A)  & 1 & 1 + 1 + 1 & (2B, 2B, 2B) & 3 \\
& & &  1 & 1 + 1 + 1 & (2B, 2A, 2A) & 6 \\
&&&&&& \\
2^3  & 4 + 2 + 2 &  ( 4B, 4B, 2B, 2A, 2A) & 1 & 1 + 1 + 1 & (2B, 2A, 2A) & 6 \\
&&&&&& \\
&&&&&& \\
S_4 & 6 + 6 & (2B, 2B, 2A)  & S_4  &  6 & (2B, 3A) & 13  \\
& &  &  1 &  1 + 1 + 1 & (2B, 2A, 2A) &  6  \\
&&&&&& \\
S_4 & 6 + 6 + 6 &  ( 2B, 2A, 2A, 2B, 2A, 2B) & 1 & 1 + 1 + 1 & (2B, 2A, 2A) & 6 \\
&&&&&& \\
S_4 & 6 + 6 + 6 &  ( 2A, 2A, 2A, 2A, 2A, 2A)  & 1 & 1 + 1 + 1 & (2A, 2A, 2A) & 9 \\
& & &  1 & 1 + 1 + 1 & (2B, 2A, 2A) & 6 \\
&&&&&& \\
S_4 & 6 + 6 + 3 &  ( 4A, 2A, 2B, 2B )  & S_4 & 6 & (2B, 3A) & 13 \\
& & &  1 & 1 + 1 + 1 & (2B, 2A, 2A) & 6 \\
& & &  1 & 1 + 1 + 1 & (2B, 2B, 2B) & 3 \\
& & &  1 & 1 + 1 + 1 & (2B, 2B, 2A) & 4 \\
&&&&&& \\
S_4 & 6 + 6 + 3 &  ( 4A, 2A, 2B, 2A)  & 1 & 1 + 1 + 1 & (2B, 2A, 2A) & 6 \\
& & &  1 & 1 + 1 + 1 & (2B, 2B, 2B) & 3 \\
&&&&&& \\
S_4 & 6 + 6 + 3 &  ( 4A, 2A, 2A, 2A)  & 1 & 1 + 1 + 1 & (2B, 2A, 2A) & 6 \\
& & &  1 & 1 + 1 + 1 & (2B, 2B, 2B) & 3 \\
&&&&&& \\ \hline
\caption{The incomplete non-minimal examples }
\label{tab:nonmin3genincomplete}
\end{longtable}

\end{document}